\documentclass[pdflatex,sn-mathphys-num]{sn-jnl}
\usepackage{graphicx}
\usepackage{multirow}
\usepackage{amsmath,amssymb,amsfonts}
\usepackage{amsthm}
\usepackage[title]{appendix}
\usepackage{xcolor}
\usepackage{textcomp}
\usepackage{manyfoot}
\usepackage{booktabs}
\usepackage{algorithm}
\usepackage{algorithmicx}
\usepackage{algpseudocode}
\usepackage{listings}

\usepackage{stmaryrd}
\DeclareFontFamily{U}{stmry}{}
\DeclareFontShape{U}{stmry}{m}{n}{<-6>stmary5<6-7>stmary6<7-8>stmary7<8-9>stmary8<9-10>stmary9<10->stmary10}{}
\DeclareFontShape{U}{stmry}{b}{n}{<->ssub * stmry/m/n}{}
\hypersetup{bookmarksdepth=subsubsection}
\theoremstyle{thmstyleone}

\newtheorem{lemma}{Lemma}
\newtheorem{proposition}{Proposition}
\newtheorem{assumption}{Assumption}
\newtheorem{corollary}{Corollary}
\newtheorem{thm}{Theorem}
\theoremstyle{thmstylethree}
\newtheorem{remark}{Remark}

\newcommand{\proofstep}[1]{\par\medskip\noindent\emph{#1}}

\newcommand{\mA}{\mathbf{A}}
\newcommand{\vb}{\mathbf{b}}
\newcommand{\vx}{\mathbf{x}}
\newcommand{\va}{\mathbf{a}}
\newcommand{\ve}{\mathbf{e}}
\newcommand{\vr}{\mathbf{r}}

\begin{document}

\title[Averaged quantile Kaczmarz methods beyond the error horizon]{Beyond the error horizon: averaged quantile Kaczmarz methods for corrupted and noisy linear systems}

\author[1]{\fnm{Sofiia} \sur{Shvaiko}}\email{ss7950@princeton.edu}
\author*[2]{\fnm{Elizaveta} \sur{Rebrova}}\email{elizaveta.rebrova@maths.ox.ac.uk}

\affil[1]{\orgdiv{ORFE}, \orgname{Princeton University}, \orgaddress{\city{Princeton}, \postcode{08544}, \country{United States}}}
\affil*[2]{\orgdiv{Mathematical Institute}, \orgname{University of Oxford}, \orgaddress{\city{Oxford}, \postcode{OX2 6GG}, \country{United Kingdom}}}

\abstract{Linear systems and least-squares problems can contain measurement
errors of different scales: a few measurements may be highly
unreliable, while small noise affects every measurement. Residual-based iterative filtering can remove the influence of sparse
corruption, but its accuracy is limited by a noise-dependent
error horizon. We study how averaging and relaxation reduce this horizon in
averaged quantile randomized Kaczmarz (AQRK), which filters equations by
residual quantiles and averages a batch of single-row updates. 
Filtering, however, introduces an update bias that averaging alone cannot remove. We show that re-querying the accepted equations can remove this bias, and  
when the corruption is sparse enough
relative to the conditioning of the matrix and large enough relative to the
noise, we prove that AQRK with re-queried residuals reduces the horizon as averaging does without the corruption: its error decays to a horizon proportional to $\sqrt{\alpha/s}$, where $\alpha$ is the stepsize and $s$ the update
batch size. Further, decreasing stepsizes or increasing batches gives convergence to the solution of the clean system. 
Numerical
experiments support the theory and examine settings beyond it,
including noise that cannot be re-queried and the corrupted low-rank tensor
fitting problem. In the latter, AQRK inner solves succeed where
standard alternating least-squares solvers and single-row robust updates
fail.}



\makeatletter
\patchcmd{\@maketitle}{\removelastskip\vskip36pt}{%
  \enlargethispage{2\baselineskip}%
  {\renewcommand\thefootnote{}\footnotetext{\fontsize{8bp}{9.5bp}\selectfont The authors were partially supported by National Science Foundation DMS-2309685.}}%
  \removelastskip\vskip36pt}{}{\PackageWarning{main}{Funding footnote patch failed}}
\makeatother

\maketitle

\section{Introduction}
\label{sec:intro}

Measurements that enter a linear system are rarely of uniform quality. A few of
them may be significantly wrong, for example because of a faulty sensor or a
corrupted transmission, while all of them carry small noise. The two kinds of
error motivate different treatment: a handful of large corruptions, which need
not be centered or independent of the coefficient matrix, can pull the
least-squares solution far from the desired one, whereas dense noise cannot be
removed from individual measurements and limits the attainable accuracy.

Row-action methods offer a natural way to treat both. The randomized Kaczmarz
method (RK) repeatedly samples an equation and projects the current iterate
onto its solution hyperplane~\cite{kaczmarz1973,strohmer2007randomized}. Its
quantile-based variants accept an equation only if its absolute residual lies
below a $q$-quantile of the residual magnitudes, which filters out large
corruptions~\cite{haddock2022quantilebased,steinerberger2021quantile}; with
dense noise, however, they converge only up to a noise-dependent \emph{error
horizon}~\cite{jarman2021quantilerk,battaglia2025qrkdoublequantilerkerror,coria2024quantile}.
Without corruption, such a horizon can be reduced by averaging several row
updates or by decreasing the
stepsize~\cite{moorman2021randomized,marshall2022optimal,tondji2023adaptive,tondji2026accelerated},
so it is natural to combine filtering with averaging; Section~\ref{sec:related}
reviews this literature in more detail.

The combination is subtler than it seems. Averaging reduces fluctuations but
not a bias, and filtering can create one: an equation is accepted because its
observed residual is small, so equations whose noise happens to shrink the
residual are accepted more often than those whose noise enlarges it. The
accepted noise therefore need not remain centered, and updates built from the
accepted observations can be biased even at the exact solution of an
uncorrupted system, as it was observed in~\cite{shen2019learning,koloskova2023revisiting}. Therefore, removing
corrupted equations and keeping the update noise centered are separate
objectives.

This paper studies averaged quantile randomized Kaczmarz (AQRK), which filters
equations by residual quantiles, averages a batch of their single-row updates,
and, if possible due to the nature of the problem, queries the selected equations again for the update,
so that the update noise is independent of the selection event. 
Correspondingly, we analyze a model in which residual measurements can be repeated: the
sparse corruptions remain fixed, while each query receives independent,
bounded, mean-zero noise. When the corruption is sparse enough relative to the
conditioning of the matrix and large enough relative to the noise,
Theorem~\ref{thm:main} gives, for every admissible stepsize $\alpha$ and $s$
independent update queries,
\begin{equation}
\label{eq:intro-main-result}
\mathbb E\|\ve_k\|^2
\le(1-\eta\alpha)^k\|\ve_0\|^2+\frac{\alpha\,r_\infty^2}{\eta\,s},
\end{equation}
where $\ve_k=\vx_k-\vx^\star$ is the error with respect to the clean solution
$\vx^\star$, $r_\infty$ bounds the query noise, and $\eta>0$ collects the
conditioning constants. The rate constant does not depend on $\alpha$ or $s$,
so the horizon is proportional to $\sqrt{\alpha/s}$, as for averaging without
corruption~\cite{moorman2021randomized,tondji2026accelerated}. The proof
combines control of the updates while corrupted rows may still be selected with
the sharper noise bound after filtering excludes them.

Least-squares problems also serve as building blocks inside larger numerical
methods. In alternating least-squares methods, for example, an inner solver
must control unreliable equations while retaining enough accuracy to support
the outer computation, and fresh observations are often unavailable. We test AQRK on fixed noisy systems, where the noise cannot be
refreshed and the bias described above may persist, and as an inner solver in
canonical polyadic (CP) tensor fitting~\cite{kolda2009tensor}, where the
subproblems change with the factors and can be  inconsistent even for an exact
low-rank target. Our main contributions are as follows.
\begin{itemize}
\item \textbf{Breaking the error horizon in the presence of sparse
corruption.} 
Theorem~\ref{thm:main} shows that averaging over $s$ fresh samples and relaxation with step size $\alpha$ makes the error
horizon of AQRK proportional to $\sqrt{\alpha/s}$ despite sparse corruption, and so it can be
made arbitrarily small. In particular, for Gaussian-like data, no relaxation and the batch size of the order of data dimension, the horizon is of order $ r_{\infty},$ compared with $\sqrt n\,r_\infty$ for
QuantileRK~\cite{jarman2021quantilerk}.

\item \textbf{Convergence to zero despite persistent corruption and query
noise.} Corollaries~\ref{cor:diminishing_stepsize}
and~\ref{cor:increasing_batch} give decreasing-stepsize and increasing-batch
schedules that drive the mean-square error to zero, which answers, in the
repeated-observation model, an open question of~\cite{coria2024quantile}.

\item \textbf{Filtered, averaged updates beyond the model.} On fixed noisy
systems, averaging can still lower the observed error floor relative to
QuantileRK in practice, although in theory it cannot remove the bias of fixed noise.
In CP tensor fitting, averaged quantile inner solves recover the tensor in
cases where single-row QuantileRK with the same stepsize and stopping rule either does
not, or incurs much higher reconstruction error (Section~\ref{subsec:cp_application}).
\end{itemize}
The rest of the paper is organized as follows: Section~\ref{sec:algorithm} introduces the problem setup and AQRK,
Section~\ref{sec:related} discusses related work, Section~\ref{sec:convergence}
contains the theoretical analysis, and Section~\ref{sec:experiments} contains the numerical
experiments.

\section{Problem setup and AQRK algorithm}
\label{sec:algorithm}

\subsection{Problem setup}
\label{subsec:problem_setup}
Let $\mA\in\mathbb R^{m\times n}$ be a known measurement matrix with rows
$\va_i^T$, and let $\vx^\star\in\mathbb R^n$ be the unknown reference target.
We assume $\|\va_i\|_2=1$, normalizing each equation if needed.

We model the observed right-hand side as $\vb=\vb^\star+\boldsymbol\rho+\vr$,
where the exact right-hand side $\vb^\star=\mA\vx^\star$ is also unknown.
The observations contain two types of error: a sparse vector
$\boldsymbol\rho\in\mathbb R^m$ of possibly large corruptions, with
$\|\boldsymbol\rho\|_0\le\beta m$, $\beta\in(0,1)$, and a dense perturbation
$\vr\in\mathbb R^m$. The latter may represent measurement error incurred each
time an observation is made, model mismatch, or even the inconsistency of a
least-squares fit. For example, $\vx^\star$ may be the least-squares solution for the data before sparse corruption, with $\vr$ the residual of that fit, as in the ALS subproblems of Section~\ref{subsec:cp_application}, or the true image in an inverse problem. Our aim is to approximate this
target from the corrupted and noisy observations.

When measurements can be repeated, we write
$b_i^{(\ell)}=b_i^\star+\rho_i+r_i^{(\ell)}$ for observation $\ell$ of row $i$;
the corruption remains fixed, while the dense perturbation may change between
observations. When they cannot, as for model mismatch or a least-squares
residual, we set $r_i^{(\ell)}=r_i$ for all $\ell$. For a current iterate
$\vx$, each call to the residual oracle uses a new observation of row $i$ and
returns
\begin{equation}
\label{eq:oracle-notation}
\mathcal D(i,\vx)=b_i^{(\ell)}-\langle\va_i,\vx\rangle
=b_i^\star+\rho_i+r_i^{(\ell)}-\langle\va_i,\vx\rangle.
\end{equation}

\subsection{Averaged quantile updates}
\label{subsec:aqrk}
We now present the averaged quantile randomized Kaczmarz (AQRK) algorithm. It
builds on randomized Kaczmarz (RK)~\cite{strohmer2007randomized}, which, for
unit-norm rows, samples a row $i$ uniformly and projects the current iterate
onto the hyperplane of the $i$-th equation,
\[
\vx_{k+1}=\vx_k+\bigl(b_i-\langle\va_i,\vx_k\rangle\bigr)\,\va_i .
\]
QuantileRK~\cite{haddock2022quantilebased} performs this update only if the
absolute residual of the sampled equation does not exceed a $q$-quantile of
the absolute residuals, and averaged RK (ARK)~\cite{moorman2021randomized}
replaces the single projection by an average of several, possibly relaxed,
projections. AQRK filters equations as QuantileRK does and averages a batch of
their updates as ARK does; in addition, it can query the selected equations
again for the update.

Starting from a current iterate $\vx_k$, the first sampling stage (Sample1)
chooses a uniformly random set $C_k\subseteq[m]$ of $m_0$ distinct candidate
rows without replacement, and computes and stores
$d_{k,i}=\mathcal D(i,\vx_k)$ for each $i\in C_k$.
For $q\in(0,1)$ and $h=\lceil qm_0\rceil$, define
\begin{equation}
\label{eq:tau-notation}
    \tau_{k,q}=\text{the $h$-th smallest value of }
        \{|d_{k,i}|:i\in C_k\}.
\end{equation}
The eligible set $L_k$ consists of the $h$ candidates with the smallest
absolute residuals. In the absence of ties,
\begin{equation}
\label{eq:admissible}
    L_k=\{i\in C_k:|d_{k,i}|\le\tau_{k,q}\};
\end{equation}
ties at the threshold are resolved at random so that $|L_k|=h$.

The second sampling stage (Sample2) draws $s_k$ update rows
$p_{k,1},\ldots,p_{k,s_k}$ independently and uniformly from $L_k$, with
replacement, and AQRK averages their residual-based directions,
\begin{equation}
\label{eq:aqrk-update}
    \vx_{k+1}=\vx_k+\frac{\alpha_k}{s_k}\sum_{j=1}^{s_k}
        \widehat d_{k,j}\,\va_{p_{k,j}}.
\end{equation}
The \emph{fresh} and \emph{reused} variants of AQRK differ only in the update
residuals: the fresh variant re-queries the oracle,
$\widehat d_{k,j}=\mathcal D(p_{k,j},\vx_k)$, whereas the reused variant takes
$\widehat d_{k,j}=d_{k,p_{k,j}}$ from the selection stage. When measurements
cannot be repeated, the two variants coincide. See Algorithm~\ref{alg:framework}.

\begin{algorithm}[H]
\caption{AQRK sampled, with or without re-querying}
\label{alg:framework}
\begin{algorithmic}[1]
\Require Matrix $\mA\in\mathbb R^{m\times n}$ with unit-norm rows $\va_i\in\mathbb R^n$, residual oracle $\mathcal D$, initial iterate $\vx_0\in\mathbb R^n$
\Require Candidate sample size $m_0\in[m]$, quantile $q\in(0,1)$, batch sizes $s_k\ge1$, stepsizes $\alpha_k>0$, iteration count $K$
\Require Re-query update residuals: \emph{fresh} or \emph{reused}
\State $h\gets\lceil qm_0\rceil$
\For{$k=0,\ldots,K-1$}
    \State Choose $C_k$ uniformly among the $m_0$-element subsets of $[m]$ \Comment{Sample1}
    \State Query and cache $d_{k,i}\gets\mathcal D(i,\vx_k)$ for every $i\in C_k$
    \State Let $L_k$ index the $h$ smallest $|d_{k,i}|$, breaking ties at random \Comment{Quantile}
    \State Draw $p_{k,1},\ldots,p_{k,s_k}\stackrel{\rm iid}{\sim}\mathrm{Uniform}(L_k)$ \Comment{Sample2}
    \For{$j=1,\ldots,s_k$}
        \If{re-querying is \emph{fresh}}
            \State $\widehat d_{k,j}\gets\mathcal D(p_{k,j},\vx_k)$ 
        \Else
            \State $\widehat d_{k,j}\gets d_{k,p_{k,j}}$
        \EndIf
    \EndFor
    \State $\displaystyle \vx_{k+1}\gets\vx_k+\frac{\alpha_k}{s_k}\sum_{j=1}^{s_k}\widehat d_{k,j}\va_{p_{k,j}}$ \Comment{Average}
\EndFor
\State \Return $\vx_K$
\end{algorithmic}
\end{algorithm}

Algorithm~\ref{alg:framework} can run for a predefined number of iterations
or until another stopping criterion is met; see, for example, the numerical
experiments in Section~\ref{sec:experiments}.

Without re-querying -- to save computational work or when fresh observations
are not available -- we can instead average all $h$ eligible rows using the
residuals already computed for selection, replacing Sample2 and the update in
Algorithm~\ref{alg:framework} by
\begin{equation}
\label{eq:aqrk-all}
    \vx_{k+1}=\vx_k+\frac{\alpha_k}{h}\sum_{i\in L_k}d_{k,i}\,\va_i .
\end{equation}
This is the QABK update of \cite{cheng2022block}, and we call it AQRK-all
here for notational consistency. Note that taking $s_k=h$ in
Algorithm~\ref{alg:framework} does not give exactly the same update, since
sampling with replacement can repeat some rows and omit others. The CP
experiments in Section~\ref{subsec:cp_application} also show that sampled
updates can be more stable than all-row averaging.

With full screening, QuantileRK is the reused variant of
Algorithm~\ref{alg:framework} with $s_k=1$ and $\alpha_k=1$, up to idle
iterations, which occur with probability $1-h/m$; RK and ARK correspond to
sampling the update rows uniformly from $[m]$ instead of $L_k$.

\begin{remark}[Update complexity and role of the sampling stages]
\label{rem:sampling_stages}
The first sampling stage (Sample1) dominates the cost of an iteration: it
requires $m_0$ inner products, that is, $O(m_0n)$ operations, or the number
of nonzeros in the candidate rows for sparse $\mA$, after which the $h$
smallest absolute residuals can be selected in $O(m_0)$ expected time.
Choosing $m_0<m$ reduces this cost; for subsampling strategies and their
guarantees, we refer to
\cite{haddock2022quantilebased,cheng2022block,haddock2023subsampled,cai2026subsample,haddock2026lossfiltering,battaglia2026streaming}.
For example, the deterministic screening conditions in
\cite{haddock2026lossfiltering} require $m_0\gtrsim n+\beta m$. Here we
focus on attainable accuracy and do not analyze candidate subsampling.
When the corruptions are sufficiently large relative to the dense noise,
corrupted observations have larger absolute residuals as the iterates
approach $\vx^\star$ and are filtered out. With full screening, $m_0=m$,
the condition $q<1-\beta$ gives $h\le m-\|\boldsymbol\rho\|_0$, so there
are enough uncorrupted equations to fill the eligible set.

The second sampling stage (Sample2) draws update rows with replacement, so
$s_k$ may exceed $h$, and repeated occurrences of a row provide independent
observations when each is queried afresh. Re-querying avoids the selection
bias described in Section~\ref{sec:intro}. It costs $s_k$ additional
observations but no new inner products, since $\langle\va_i,\vx_k\rangle$ can
be cached during selection, and grouping repeated rows keeps the update cost
at $O(\min\{s_k,h\}\,n+s_k)$ operations. Thus, with full screening, an iteration costs $m+s_k$ residual queries and
$O(mn+s_k)$ operations. We return to the cost of the batch size in
Section~\ref{subsec:cost}, where we compare decreasing stepsizes with
increasing batches. 
\end{remark}

\section{Related work}
\label{sec:related}
Randomized Kaczmarz (RK)~\cite{kaczmarz1973,strohmer2007randomized}
converges exponentially in expectation on consistent systems; with a
noisy right-hand side, it converges only up to a horizon
$\sqrt{R}\,\gamma$~\cite{needell2010}, where
$R=\|{\bf A}\|_F^2/\sigma_{\min}({\bf A})^2\ge n$ and $\gamma$ is the
largest row-scaled noise entry. When the noise is centered and
independent of the past iterates, this horizon can be reduced or removed in
two ways.

The first is to decrease the stepsizes, which recover the clean solution to
any accuracy within a single pass without replacement, given enough
equations~\cite{marshall2022optimal}, or, when every query returns fresh noise,
over arbitrarily many iterations~\cite{tondji2023adaptive}. The second is
averaging: averaging $\tau$ row updates reduces the noise variance by a factor
of order $1/\tau$~\cite{moorman2021randomized}, and averaged adaptive
Bregman--Kaczmarz combines averaging over fresh, draw-indexed queries with
adaptive stepsizes to obtain convergence to zero at a rate that improves with
the batch size~\cite{tondji2026accelerated}; see
also~\cite{lin2015learning,2025tailaveraging,zhang2026signal}. These methods
use every sampled equation and do not treat sparse corruption.

For a right-hand side with a fraction $\beta$ of arbitrarily corrupted
entries, QuantileRK~\cite{haddock2022quantilebased} updates only with
equations whose absolute residual lies below the $q$-quantile of the
residuals, and for suitable random matrices it satisfies
$\mathbb E\|\ve_k\|^2\le(1-C_q/n)^k\|\ve_0\|^2$, where
$\ve_k=\vx_k-\vx^\star$. Steinerberger~\cite{steinerberger2021quantile}
gives a deterministic guarantee under a lower bound on
$\sigma_{\min,q-\beta}^2/\sigma_{\max}^2$, of the same type as our
condition $\mu>0$, and Haddock and Needell~\cite{haddock2019randomized} show
that near the clean solution the largest residuals identify the corrupted
equations. Our direct predecessor is the quantile averaged block Kaczmarz
method of Cheng et al.~\cite{cheng2022block}, AQRK-all in
Section~\ref{subsec:aqrk}. Its analysis assumes exact clean equations, so no
noise horizon arises. Further work reduces the cost of the quantile by
subsampling \cite{haddock2023subsampled}, removes persistently
unreliable rows using martingale techniques from~\cite{anderson2025beyond}
to show that the largest residuals identify the corruptions after a
warm-up~\cite{shvaiko2026whitelist}, and extends quantile filtering to tensor
systems~\cite{castillo2025tensor} and nonnegative matrix
factorization~\cite{collins2025quantile}.

Closest to our setting are analyses of QuantileRK with both sparse corruption
and dense bounded noise $\vr$. For random matrices, Jarman and
Needell~\cite{jarman2021quantilerk} prove
\[
\mathbb E\|\ve_k\|^2\le\Bigl(1-\frac{Cq}{n}\Bigr)^k\|\ve_0\|^2
+\frac{2n}{C'}\,\|\vr\|_\infty^2 ,
\]
and Battaglia and Ma~\cite{battaglia2025qrkdoublequantilerkerror} obtain the
deterministic bound
\[
\mathbb E\|\ve_k\|^2\le(1-C)^k\|\ve_0\|^2
+\frac1C\Bigl(\frac{2\beta(1-q)}{q(1-q-\beta)}+1\Bigr)\|\vr\|_\infty^2 ,
\]
with $C$ depending on $q$, $\beta$ and the singular values of $\mA$, together
with bounds for a double-quantile variant. Coria, Haddock, and
Pacheco~\cite{coria2024quantile} allow both the noise and the corruption to
change between iterations, and Zhang, Wang, and
Zhang~\cite{zhang2025quantilesparse} combine quantile filtering with averaged
sparse Kaczmarz updates, later accelerated by adaptive
momentum~\cite{zhang2026momentum}. In all of these analyses, the filter and
the update use the same observation, and the bounds keep a horizon that
neither smaller stepsizes nor larger batches remove: the single-row analyses
use unit steps, and the averaged bound of~\cite{zhang2025quantilesparse}
retains a noise term linear in the stepsize. Coria et al.\ pose as an open
question a Kaczmarz-type method that is robust to corruption and noise
simultaneously and breaks this horizon.
We answer it under repeated
observations: with fresh update queries, the squared horizon of AQRK
(as in~\eqref{eq:intro-main-result}) is proportional to $\alpha/s$, and
decreasing-stepsize and increasing-batch schedules drive the mean-square
error to zero. The price is a lower bound on the corruption magnitudes, which
the bounds above do not require. 

Instead of filtering equations, one can replace the least-squares objective
by a robust loss that reduces the influence of large residuals; Schöpfer et
al.~\cite{schopfer2022extended} minimize such losses for sparse large errors
with an extended randomized Kaczmarz method. Their guarantees concern the
minimizer of the robust objective for the observed data, rather than the
solution of the clean system. Robustness can also come from controlling step
lengths instead of rejecting observations. QuantileSGD~\cite{haddock2022quantilebased}
takes sign-based $\ell_1$ steps whose length is set by a residual quantile,
Jeong, Needell, and Rebrova~\cite{jeong2024stochastic} use exponentially
decreasing step lengths for streaming systems with Massart corruption, and
Pesme and Flammarion~\cite{pesme2020online} analyze $\ell_1$-loss SGD with
decreasing stepsizes and iterate averaging when the sparse corruptions are
independent of Gaussian covariates and dense noise. In our approach, the filter
provides the robustness to corruption, and decreasing stepsizes serve a
different purpose: combined with filtering and fresh update queries, they drive
the mean-square error to zero
(Corollary~\ref{cor:diminishing_stepsize}).

\section{Convergence of AQRK}
\label{sec:convergence}

\subsection{Setup for the analysis}
\label{subsec:analysis_setup}
Throughout this section, we analyze Algorithm~\ref{alg:framework} with
re-querying and full screening, $C_k=[m]$ at every iteration. We write $\ve_k=\vx_k-\vx^\star$ for the error.

\begin{assumption}[Noise]
\label{as:model}
The perturbations $r_i^{(\ell)}$ returned by different oracle calls are
independent of each other and of the sampling in
Algorithm~\ref{alg:framework}, have mean zero, and satisfy
$|r_i^{(\ell)}|\le r_\infty$ almost surely.
\end{assumption}

At iteration $k$, we write $r^{\rm qnt}_{k,i}$ for the noise in the selection
observation of row $i$ and $r^{\rm upd}_{k,j}$ for the noise in the $j$-th
update observation, so that
\[
d_{k,i}=b_i^\star+\rho_i-\langle\va_i,\vx_k\rangle+r^{\rm qnt}_{k,i},
\qquad
\widehat d_{k,j}=b_{p_{k,j}}^\star+\rho_{p_{k,j}}
-\langle\va_{p_{k,j}},\vx_k\rangle+r^{\rm upd}_{k,j}.
\]
Note that a row drawn
twice during the update phase also receives two independent noise values.

\begin{assumption}[Corruption and quantile]
\label{as:high_separation}
Every nonzero corruption satisfies $|\rho_i|\ge C_{\rm corr}$. The quantile
level satisfies $\beta<q<1-\beta$, and $qm$ and $\beta m$ are integers.
\end{assumption}

We write $\sigma_{\max}=\sigma_{\max}(\mA)$ and
\[
\sigma_{\min,q-\beta}=\min_{T\subseteq[m],\ |T|=(q-\beta)m}\sigma_{\min}(\mA_T),
\]
where $\sigma_{\min}(\mA_T)^2=\lambda_{\min}(\mA_T^T\mA_T)$, which is zero when
$\mA_T$ does not have full column rank.

 \subsection{Main result}
\label{subsec:main_result}
We define
\[
\omega=\frac{1}{\sqrt{1-q-\beta}},\qquad
\kappa=\frac{\sigma_{\max}^2}{qm},\qquad
\mu=\sigma_{\min,q-\beta}^2-2\sqrt\beta\,\omega\,\sigma_{\max}^2 .
\]
Intuitively, $\mu$ is the smallest curvature over subsets of $(q-\beta)m$
rows, reduced by a penalty for corrupted rows that pass the filter. For a
margin $\gamma\in(0,1/2]$ and a balancing parameter $\delta>1$, we define
\begin{equation}
\label{eq:main_constants}
\begin{aligned}
\eta&=\frac{2(1-\gamma)(1-\delta^{-1})\,\mu}{qm},\\
\alpha_{\max}(s)&=\frac{2\gamma\,\mu/\sigma_{\max}^2}
{\dfrac{1+4\beta\omega^2}{s}+\kappa
+\dfrac{4\beta\omega(\beta\omega+\sqrt{q\kappa})}{q}},\\
\xi_s^2&=\frac{8\delta}{\omega(1-\gamma)}\,\frac{m\,r_\infty^2}{\mu}
\left(\sqrt\beta+\beta\,\alpha_{\max}(s)
\left(\frac{\omega}{s}+\frac{\beta\omega+\sqrt{q\kappa}}{q}\right)\right).
\end{aligned}
\end{equation}
Here $\eta$ is the rate constant, $\alpha_{\max}(s)$ bounds the stepsize, and
the radius $\xi_s$ determines how large the corruptions must be compared with
the noise for the filtering to be provably effective. Our central
convergence result is as follows.

\begin{thm}
\label{thm:main}
Let Assumptions~\ref{as:model} and~\ref{as:high_separation} hold, and suppose
that $\mu>0$. Fix $\gamma\in(0,1/2]$, $\delta>1$, $s\ge1$ and
$0<\alpha\le\alpha_{\max}(s)$, and run Algorithm~\ref{alg:framework} with
$s_k\equiv s$ and $\alpha_k\equiv\alpha$. If the separation condition
\begin{equation}
\label{eq:separation}
C_{\rm corr}>2(r_\infty+\xi_s)
\end{equation}
holds, then $0<\eta\alpha<1$ and, for all $k\ge0$,
\[
\mathbb E\|\ve_k\|^2
\le(1-\eta\alpha)^k\|\ve_0\|^2+\frac{\alpha\,r_\infty^2}{\eta\,s}.
\]
\end{thm}

The bound consists of a geometric term and a horizon. The rate constant
$\eta$ depends on $\mA$, $q$ and $\beta$ through $\mu/(qm)$, but not on $s$,
$\alpha$ or the noise level. For a given stepsize, the (squared) horizon
$\alpha r_\infty^2/(\eta s)$ is proportional to $\alpha/s$: averaging over
more fresh observations and decreasing the stepsize reduce it in the same
way, while a smaller stepsize also slows the geometric decay and a larger $s$
requires more fresh observations per iteration. Since $\alpha_{\max}(s)$ is
nondecreasing in $s$, any stepsize admissible for single-row updates,
$\alpha\le\alpha_{\max}(1)$, remains admissible for every batch size.

The condition $\mu>0$ limits the
corruption fraction relative to the restricted conditioning of $\mA$. The parameters $\gamma$ and
$\delta$ tune the constants. The choice $\gamma=1/2$ admits the largest
steps, while a smaller $\gamma$ restricts the stepsize but improves $\eta$
and $\xi_s$ by constant factors. Increasing $\delta$ improves $\eta$ but
enlarges $\xi_s$, and thus balances the rate against the required
separation.

The proof in Section~\ref{subsec:proofs} analyzes two regimes,
separated by the radius $\xi_s$. Far from $\vx^\star$, corrupted rows may
pass the filter, but their contribution is controlled through the quantile
threshold and absorbed by a fraction $1/\delta$ of the contraction whenever
$\|\ve_k\|>\xi_s$. Near $\vx^\star$, the separation
condition~\eqref{eq:separation} ensures that every corrupted residual exceeds
every uncorrupted one, so the filter rejects all corrupted rows and only the
fresh noise remains. Both regimes lead to the one-step recursion of
Proposition~\ref{prop:onestep}, which we iterate.

\subsection{Proof of Theorem~\ref{thm:main}}
\label{subsec:proofs}
We first present two auxiliary lemmas based on known techniques: the first
bounds the quantile threshold, and the second shows that corrupted rows are
filtered out near $\vx^\star$. We then bound one step of AQRK and combine
these estimates into a one-step recursion. 

Lemma~\ref{lem:quantile} bounds the threshold $\tau_{k,q}$. For
$r_\infty=0$, it is~\cite[Lemma~1]{steinerberger2021quantile}. Under our
assumptions, it sharpens the noisy bounds
of~\cite[Corollary~II.4]{jarman2021quantilerk}, which hold with high
probability for random matrices, and~\cite[Lemma~2.2]{coria2024quantile}.
\begin{lemma}[Quantile bound]
 \label{lem:quantile}
At iteration $k$, recall that we denote as $\tau_{k,q}$ the $qm$-th smallest in absolute value residual from the selection set $|d_{k,1}|,\ldots,|d_{k,m}|$, see~\eqref{eq:tau-notation}.
Then, with $\omega=1/\sqrt{1-q-\beta}$,
\[
\tau_{k,q}\le\frac{\omega\,\sigma_{\max}\|\ve_k\|}{\sqrt m}+r_\infty .
\]
\end{lemma}

\begin{proof}
Let $U$ denote the set of
uncorrupted rows. By definition of $\tau_{k,q}$, at
least $m-qm+1\ge(1-q)m$ rows satisfy $|d_{k,i}|\ge\tau_{k,q}$. At most
$\beta m$ of these rows are corrupted, so the set $G\subseteq U$ of
uncorrupted rows with $|d_{k,i}|\ge\tau_{k,q}$ satisfies
$|G|\ge(1-q-\beta)m=m/\omega^2>0$. If $\tau_{k,q}\le r_\infty$, the claim
holds. Otherwise, for $i\in G$ we have $\rho_i=0$ and
$d_{k,i}=-\langle\va_i,\ve_k\rangle+r^{\rm qnt}_{k,i}$, so
\[
\tau_{k,q}\le|d_{k,i}|\le|\langle\va_i,\ve_k\rangle|+r_\infty,
\qquad\text{that is,}\qquad
0<\tau_{k,q}-r_\infty\le|\langle\va_i,\ve_k\rangle| .
\]
Squaring and averaging over $G$ gives
\[
\begin{aligned}
(\tau_{k,q}-r_\infty)^2
&\le\frac{1}{|G|}\sum_{i\in G}\langle\va_i,\ve_k\rangle^2
\le\frac{\omega^2}{m}\sum_{i\in U}\langle\va_i,\ve_k\rangle^2\\
&=\frac{\omega^2}{m}\|\mA_U\ve_k\|^2
\le\frac{\omega^2\sigma_{\max}^2}{m}\|\ve_k\|^2 ,
\end{aligned}
\]
and taking square roots proves the claim.
\end{proof}

The next lemma extends~\cite[Lemma~1]{haddock2019randomized} to the noisy
case; see also~\cite[Lemma~3.2]{shvaiko2026whitelist}.

\begin{lemma}[Filtering]
\label{lem:filtering}
At iteration $k$, if $\|\ve_k\|<C_{\rm corr}/2-r_\infty$, then $L_k\subseteq U$, where $U$ denotes the set of all uncorrupted rows, and $L_k$ is the eligible set.
\end{lemma}

\begin{proof}
For $i\in U$, since $\|\va_i\|=1$ and $|r^{\rm qnt}_{k,i}|\le r_\infty$,
\[
|d_{k,i}|=|-\langle\va_i,\ve_k\rangle+r^{\rm qnt}_{k,i}|
\le\|\ve_k\|+r_\infty .
\]
For a corrupted row $i$, since $|\rho_i|\ge C_{\rm corr}$,
\[
|d_{k,i}|=|\rho_i-\langle\va_i,\ve_k\rangle+r^{\rm qnt}_{k,i}|
\ge C_{\rm corr}-\|\ve_k\|-r_\infty .
\]
From the assumption $\|\ve_k\|<C_{\rm corr}/2-r_\infty$ we get that every corrupted
absolute residual is larger than every uncorrupted one. Since $L_k$ consists
of the $qm$ smallest absolute residuals and $|U|\ge(1-\beta)m>qm$, it
contains no corrupted row.
\end{proof}

Next, we will state and prove the expected error improvement after one step of AQRK. 
We use the abbreviations
\begin{equation}
\label{eq:proof_constants}
\begin{aligned}
c_1&=\frac{2\mu}{qm}, &
c_2(s)&=\kappa\left(\frac{1+4\beta\omega^2}{s}+\kappa
+\frac{4\beta\omega(\beta\omega+\sqrt{q\kappa})}{q}\right),\\
c_3&=\frac{16\sqrt\beta}{q\,\omega}, &
c_4(s)&=\frac{16\beta}{q\,\omega}\left(\frac{\omega}{s}
+\frac{\beta\omega+\sqrt{q\kappa}}{q}\right).
\end{aligned}
\end{equation}
We denote by $\mathcal F_k$ the sigma-algebra generated by the initial iterate, all sampling and
observations before iteration $k$, and the selection observations at
iteration $k$, so that $\vx_k$ and $L_k$ are fixed under
$\mathbb E_k[\cdot]=\mathbb E[\cdot\mid\mathcal F_k]$. 

\begin{lemma}[One-step bound]
\label{lem:onestep_update}
At iteration $k$, let the batch size be $s$ and the stepsize $\alpha>0$.
Then
\begin{equation}
\label{eq:onestep_general}
\mathbb E_k\|\ve_{k+1}\|^2
\le\bigl(1-c_1\alpha+c_2(s)\alpha^2\bigr)\|\ve_k\|^2
+r_\infty^2\bigl(c_3\alpha+c_4(s)\alpha^2\bigr)
+\frac{\alpha^2r_\infty^2}{s}.
\end{equation}
If $L_k\subseteq U$, the set
of all uncorrupted rows, then the term with $c_3$ and $c_4(s)$ can be omitted.
\end{lemma}

\begin{proof}
Write $\ve=\ve_k$ and $L=L_k$ for the eligible set, $p_j=p_{k,j}$ for the
update rows drawn from $L$, and $r_j=r^{\rm upd}_{k,j}$ for their fresh noise.

\proofstep{Step 1: error recursion and update noise.} With
stepsize $\alpha$ and batch size $s$, the update~\eqref{eq:aqrk-update} with
fresh residuals reads
\[
\vx_{k+1}=\vx_k+\frac{\alpha}{s}\sum_{j=1}^s\widehat d_{k,j}\,\va_{p_j}.
\]
Since $b_i^\star=\langle\va_i,\vx^\star\rangle$, each update residual is
\[
\widehat d_{k,j}=b_{p_j}^\star+\rho_{p_j}-\langle\va_{p_j},\vx_k\rangle
+r_j
=-\bigl(\langle\va_{p_j},\ve\rangle-\rho_{p_j}\bigr)+r_j,
\]
so the sampled row $p_j$ moves the iterate along $-v_{p_j}+r_j\va_{p_j}$,
where
\[
v_i=(\langle\va_i,\ve\rangle-\rho_i)\,\va_i,\qquad i\in L,
\]
is the noiseless update direction of row $i$. So,
\[
\ve_{k+1}=\ve_k-\alpha\bar v+\alpha\bar r,
\qquad
\bar v=\frac1s\sum_{j=1}^s v_{p_j},
\qquad
\bar r=\frac1s\sum_{j=1}^s r_j\va_{p_j}.
\]

Conditional on $\mathcal F_k$ and $p_1,\ldots,p_s$, the
vector $\ve-\alpha\bar v$ is fixed, while, by Assumption~\ref{as:model}, the
update noises $r_1,\ldots,r_s$ are independent,
have mean zero and are bounded by $r_\infty$. Writing
$\mathbb E_{k,p}[\cdot]=\mathbb E[\cdot\mid\mathcal F_k,p_1,\ldots,p_s]$, we
obtain with $\|\va_{p_j}\|=1$,
\[
\mathbb E_{k,p}\langle\ve-\alpha\bar v,\bar r\rangle=0,
\qquad
\mathbb E_{k,p}\|\bar r\|^2
=\frac{1}{s^2}\sum_{j=1}^s\mathbb E_{k,p}\bigl[r_j^2\bigr]
\,\|\va_{p_j}\|^2
\le\frac{r_\infty^2}{s}.
\] Therefore,
\begin{equation}
\label{eq:noise_split}
\mathbb E_k\|\ve_{k+1}\|^2
\le\|\ve\|^2-2\alpha\langle\ve,\mathbb E_k\bar v\rangle
+\alpha^2\,\mathbb E_k\|\bar v\|^2+\frac{\alpha^2r_\infty^2}{s}.
\end{equation}

\proofstep{Step 2: averaging over the sampled rows.} It remains to
bound the two terms of~\eqref{eq:noise_split} that involve $\bar v$: the cross
term, which must produce the decrease, and the second moment, which must be of
order $\alpha^2$. 

Conditional on $\mathcal F_k$, the indices $p_j$ are
independent and uniform on $L$, and $|L|=qm$. With \[g=\sum_{i\in L}v_i, 
\;\; \text{ we have for each $j$ } \;\; \mathbb E_kv_{p_j}=\frac{g}{qm} \;\; \text{ and so }  \;\; \mathbb E_k \bar v = \frac{g}{qm}.\]
Note that the vector
$g/(qm)$ is the expected noiseless update direction. For
the second moment, independence of $p_j$ and $p_l$ for $j\ne l$ gives
\[
\mathbb E_k\|\bar v\|^2
=\frac1{s^2}\sum_{j=1}^s\mathbb E_k\|v_{p_j}\|^2
+\frac1{s^2}\sum_{j\ne l}\bigl\langle\mathbb E_kv_{p_j},
\mathbb E_kv_{p_l}\bigr\rangle
= \frac{1}{sqm}\sum_{i\in L}\|v_i\|^2+\frac{s-1}{s}\,\frac{\|g\|^2}{q^2m^2}.
\]
So, \eqref{eq:noise_split} is further bounded as
\begin{equation}
\label{eq:three_sums}
\mathbb E_k\|\ve_{k+1}\|^2
\le\|\ve\|^2-\frac{2\alpha}{qm}\,\langle\ve,g\rangle
+\frac{\alpha^2}{sqm}\sum_{i\in L}\|v_i\|^2
+\frac{\alpha^2}{q^2m^2}\,\|g\|^2+\frac{\alpha^2r_\infty^2}{s}.
\end{equation}
Note that the right-hand side no longer involves the sampled rows, only the eligible
set.

\proofstep{Step 3: clean and corrupted rows.} The decrease comes
from $\langle\ve,g\rangle$, which we bound from below, while
$\sum_{i\in L}\|v_i\|^2$ and $\|g\|^2$ must be bounded from above. Each of
these is a sum over $L=L_1\cup L_2$, where we split the eligible
set into its clean rows $L_1=L\cap U$ and its corrupted rows
$L_2=L\setminus U$. The two parts behave differently:
over clean rows, the sums are quadratic forms in $\ve$, whereas over corrupted
rows they are controlled only through the eligibility threshold.

For the clean rows, let $\mathbf M=\mA_{L_1}^T\mA_{L_1}$. For $i\in L_1$, we
have $\rho_i=0$ and $v_i=\va_i\va_i^T\ve$, and the rows have unit norm, so
the three clean parts are
\[
\sum_{i\in L_1}\langle\ve,v_i\rangle=\ve^T\mathbf M\ve,
\qquad
\sum_{i\in L_1}\|v_i\|^2=\sum_{i\in L_1}\langle\va_i,\ve\rangle^2
=\ve^T\mathbf M\ve,
\qquad
\sum_{i\in L_1}v_i=\mathbf M\ve .
\]

For the corrupted rows $i\in L_2$, eligibility gives
$|d_{k,i}|\le\tau_{k,q}$ and so
\[
|\langle\va_i,\ve\rangle-\rho_i|=|-d_{k,i}+r^{\rm qnt}_{k,i}|
\le\tau_{k,q}+r_\infty=:Q,
\]
where $Q$ bounds the length of the update direction of every corrupted row
that passes the filter.
Hence $\|v_i\|\le Q$ for $i\in L_2$, and since $|L_2|\le\beta m$, the three
corrupted parts satisfy
\[
\sum_{i\in L_2}\|v_i\|^2\le\beta m\,Q^2,
\qquad
\Bigl\|\sum_{i\in L_2}v_i\Bigr\|\le\beta m\,Q,
\]
\[
\begin{aligned}
\Bigl|\sum_{i\in L_2}\langle\ve,v_i\rangle\Bigr|
&=\Bigl|\sum_{i\in L_2}(\langle\va_i,\ve\rangle-\rho_i)
\langle\va_i,\ve\rangle\Bigr|
\le\sqrt{\beta m}\,Q\,\sigma_{\max}\|\ve\| 
\end{aligned}
\]
by the Cauchy--Schwarz inequality and since $\mA_{L_2}$ is a row
submatrix of $\mA$. Altogether, 
\begin{equation}
\label{eq:combined_bound}
\begin{aligned}
\mathbb E_k\|\ve_{k+1}\|^2
\le{}&\|\ve\|^2-\frac{2\alpha}{qm}\,\ve^T\mathbf M\ve
+\frac{2\alpha\sqrt{\beta m}\,\sigma_{\max}}{qm}\,\|\ve\|\,Q\\
&+\frac{\alpha^2}{sqm}\bigl(\ve^T\mathbf M\ve+\beta m\,Q^2\bigr)\\
&+\frac{\alpha^2}{q^2m^2}\bigl(\|\mathbf M\ve\|^2
+2\beta m\,Q\,\|\mathbf M\ve\|+\beta^2m^2Q^2\bigr)
+\frac{\alpha^2r_\infty^2}{s}.
\end{aligned}
\end{equation}

\proofstep{Step 4: matrix bounds.} The bound~\eqref{eq:combined_bound}
still depends on the eligible set through $\mathbf M$. We replace it by
bounds that hold for every eligible set, so that only the threshold~$Q$
remains to be controlled. 

Since $\mA_{L_1}$ is a row submatrix
of $\mA$, $\ve^T\mathbf M\ve\le\sigma_{\max}^2\|\ve\|^2$ and
$\|\mathbf M\ve\|\le\sigma_{\max}^2\|\ve\|$. For the lower bound, note that
$|L_1|=|L|-|L_2|\ge(q-\beta)m$ and that removing rows can only decrease
$\ve^T\mathbf M\ve=\|\mA_{L_1}\ve\|^2$. Applying the definition of
$\sigma_{\min,q-\beta}$ to any $(q-\beta)m$ of the rows in $L_1$ therefore
gives $\ve^T\mathbf M\ve\ge\sigma_{\min,q-\beta}^2\|\ve\|^2$.  With
$\kappa = \sigma_{\max}^2/(qm)$,~\eqref{eq:combined_bound} simplifies to
\begin{equation}
\label{eq:pre_quantile}
\begin{aligned}
\mathbb E_k\|\ve_{k+1}\|^2
&\le\Bigl(1-\frac{2\alpha\sigma_{\min,q-\beta}^2}{qm}
+\alpha^2\kappa\Bigl(\frac1s+\kappa\Bigr)\Bigr)\|\ve\|^2
+\frac{\alpha^2\beta}{q}\Bigl(\frac1s+\frac\beta q\Bigr)Q^2\\
&\qquad+\frac{2\alpha}{q}\bigl(\sqrt{\beta q\kappa}+\alpha\beta\kappa\bigr)
\|\ve\|\,Q+\frac{\alpha^2r_\infty^2}{s}.
\end{aligned}
\end{equation}

If $L_k\subseteq U$, then $L_2=\varnothing$, so the corrupted parts in
Step~3 vanish and~\eqref{eq:combined_bound} holds without the terms
containing $Q$. The same substitutions then give
\[
\mathbb E_k\|\ve_{k+1}\|^2
\le\Bigl(1-\frac{2\alpha\sigma_{\min,q-\beta}^2}{qm}
+\alpha^2\kappa\Bigl(\frac1s+\kappa\Bigr)\Bigr)\|\ve\|^2
+\frac{\alpha^2r_\infty^2}{s}.
\]
Since $\mu\le\sigma_{\min,q-\beta}^2$, we have
$c_1=2\mu/(qm)\le2\sigma_{\min,q-\beta}^2/(qm)$, and since all terms of
$c_2(s)$ are nonnegative, $\kappa(1/s+\kappa)\le c_2(s)$. This proves the
second claim.

\proofstep{Step 5: bound on $Q$.} It remains to control $Q$. In
general, Lemma~\ref{lem:quantile} gives
\[
Q=\tau_{k,q}+r_\infty
\le\frac{\omega\sigma_{\max}}{\sqrt m}\|\ve\|+2r_\infty
=\omega\sqrt{q\kappa}\,\|\ve\|+2r_\infty .
\]
We distinguish two cases. First, suppose that
$\omega\sqrt{q\kappa}\,\|\ve\|\ge2r_\infty$. Then
$Q\le2\omega\sqrt{q\kappa}\,\|\ve\|$, and substituting this
into~\eqref{eq:pre_quantile} gives
\[
\begin{aligned}
\mathbb E_k\|\ve_{k+1}\|^2
\le{}&\Bigl(1-\alpha\Bigl(\frac{2\sigma_{\min,q-\beta}^2}{qm}
-4\sqrt\beta\,\omega\kappa\Bigr)\\
&\quad+\alpha^2\kappa\Bigl(\frac{1+4\beta\omega^2}{s}+\kappa
+\frac{4\beta^2\omega^2}{q}+\frac{4\beta\omega\sqrt{q\kappa}}{q}\Bigr)\Bigr)
\|\ve\|^2+\frac{\alpha^2r_\infty^2}{s}.
\end{aligned}
\]
Since $\sigma_{\max}^2=qm\kappa$, the coefficient of $\alpha$ is
$2(\sigma_{\min,q-\beta}^2-2\sqrt\beta\,\omega\sigma_{\max}^2)/(qm)=c_1$, and
the coefficient of $\alpha^2$ is $c_2(s)$, which
implies~\eqref{eq:onestep_general}.

Now suppose that $\omega\sqrt{q\kappa}\,\|\ve\|<2r_\infty$. Then
$Q<4r_\infty$ and $\|\ve\|\,Q<8r_\infty^2/(\omega\sqrt{q\kappa})$, so the
terms containing $Q$ in~\eqref{eq:pre_quantile} are at most
\[
r_\infty^2\left(\frac{16\sqrt\beta}{q\,\omega}\,\alpha
+\frac{16\beta}{q\,\omega}\Bigl(\frac{\omega}{s}
+\frac{\beta\omega+\sqrt{q\kappa}}{q}\Bigr)\alpha^2\right)
=r_\infty^2\bigl(c_3\alpha+c_4(s)\alpha^2\bigr).
\]
The coefficient of $\|\ve\|^2$ in~\eqref{eq:pre_quantile} is at most
$1-c_1\alpha+c_2(s)\alpha^2$, as in the clean case. Hence
\eqref{eq:onestep_general} holds in this case as well.
\end{proof}

\begin{proposition}[One-step recursion]
\label{prop:onestep}
Let Assumptions~\ref{as:model} and~\ref{as:high_separation} hold, suppose
that $\mu>0$, and fix $\gamma\in(0,1/2]$, $\delta>1$ and $s_0\ge1$ such that
$C_{\rm corr}>2(r_\infty+\xi_{s_0})$. If, at iteration $k$, the batch size
satisfies $s_k\ge s_0$ and the stepsize satisfies
$0<\alpha_k\le\alpha_{\max}(s_0)$, then
\[
\mathbb E_k\|\ve_{k+1}\|^2
\le(1-\eta\alpha_k)\|\ve_k\|^2+\frac{\alpha_k^2r_\infty^2}{s_k}.
\]
\end{proposition}

\begin{proof} Note that $c_2(s), c_4(s)$ are both nonincreasing in $s$ and
$\alpha_k\le\alpha_{\max}(s_0)=\gamma c_1/c_2(s_0)$,
therefore
\begin{equation}
\label{eq:stepsize_effect}
1-c_1\alpha_k+c_2(s_k)\alpha_k^2\le1-(1-\gamma)c_1\alpha_k 
\end{equation}
and
\begin{equation}
\label{eq:corruption_effect}
\begin{aligned}
r_\infty^2\bigl(c_3\alpha_k+c_4(s_k)\alpha_k^2\bigr)
\le\alpha_kr_\infty^2\bigl(c_3+\alpha_{\max}(s_0)c_4(s_0)\bigr)
=\frac{(1-\gamma)c_1\alpha_k}{\delta}\,\xi_{s_0}^2 .
\end{aligned}
\end{equation}

\noindent\emph{Case $\|\ve_k\|\le\xi_{s_0}$.} The separation condition is equivalent
to $\xi_{s_0}<C_{\rm corr}/2-r_\infty$, so
$\|\ve_k\|<C_{\rm corr}/2-r_\infty$ and $L_k\subseteq U$ by
Lemma~\ref{lem:filtering}. The second part of
Lemma~\ref{lem:onestep_update} and~\eqref{eq:stepsize_effect} give
\[
\mathbb E_k\|\ve_{k+1}\|^2
\le\bigl(1-(1-\gamma)c_1\alpha_k\bigr)\|\ve_k\|^2
+\frac{\alpha_k^2r_\infty^2}{s_k}
\le(1-\eta\alpha_k)\|\ve_k\|^2+\frac{\alpha_k^2r_\infty^2}{s_k}.
\]

\noindent\emph{Case $\|\ve_k\|>\xi_{s_0}$.} The first part of
Lemma~\ref{lem:onestep_update}, \eqref{eq:stepsize_effect},
\eqref{eq:corruption_effect} and $\xi_{s_0}^2<\|\ve_k\|^2$ give
\[
\begin{aligned}
\mathbb E_k\|\ve_{k+1}\|^2
&\le\bigl(1-(1-\gamma)c_1\alpha_k\bigr)\|\ve_k\|^2
+\frac{(1-\gamma)c_1\alpha_k}{\delta}\,\xi_{s_0}^2
+\frac{\alpha_k^2r_\infty^2}{s_k}\\
&\le\Bigl(1-(1-\gamma)c_1\alpha_k+\frac{(1-\gamma)c_1\alpha_k}{\delta}\Bigr)
\|\ve_k\|^2+\frac{\alpha_k^2r_\infty^2}{s_k}\\
&=(1-\eta\alpha_k)\|\ve_k\|^2+\frac{\alpha_k^2r_\infty^2}{s_k}.
\end{aligned}
\]
\end{proof}

\begin{proof}[Proof of Theorem~\ref{thm:main}]
We apply Proposition~\ref{prop:onestep} with $s_0=s_k=s$ and
$\alpha_k=\alpha$ for all $k$. Taking total expectations gives
\begin{equation}
\label{eq:main_recursion}
\mathbb E\|\ve_{k+1}\|^2\le(1-\eta\alpha)\,\mathbb E\|\ve_k\|^2
+\frac{\alpha^2r_\infty^2}{s}.
\end{equation}
Since $\mu>0$, we have $\eta\alpha>0$. For the upper bound, $\mu\le
\sigma_{\min,q-\beta}^2\le\sigma_{\max}^2$ gives
$c_1=2\mu/(qm)\le2\sigma_{\max}^2/(qm)=2\kappa$, and $c_2(s)\ge\kappa^2$.
Hence
\[
\eta\alpha\le\eta\,\alpha_{\max}(s)
=(1-\gamma)(1-\delta^{-1})\,c_1\cdot\frac{\gamma c_1}{c_2(s)}
\le4\gamma(1-\gamma)(1-\delta^{-1})\le1-\delta^{-1}<1.
\] Iterating~\eqref{eq:main_recursion} from
$k=0$ gives
\[
\mathbb E\|\ve_k\|^2
\le(1-\eta\alpha)^k\|\ve_0\|^2
+\frac{\alpha^2r_\infty^2}{s}\sum_{i=0}^{k-1}(1-\eta\alpha)^i,
\]
and, since $0<1-\eta\alpha<1$,
$\sum_{i=0}^{k-1}(1-\eta\alpha)^i\le1/(\eta\alpha)$. This proves the
theorem.
\end{proof}

\subsection{Interpretation under random design}
\label{subsec:random_design}
The constants in Theorem~\ref{thm:main} involve restricted singular values,
which are hard to compute for a given matrix. For random matrices, they can be
estimated in terms of the dimensions, using the framework
of~\cite{haddock2022quantilebased}.

\begin{assumption}
\label{as:subgaus}
The rows $\va_i$ of $\mA$ are independent and have unit norm, and
$\sqrt n\,\va_i$ is mean zero and isotropic with subgaussian norm
$\|\sqrt n\,\va_i\|_{\psi_2}\le K$; see~\cite[Sections~2.5
and~3.4]{vershynin2018high} for the definitions.
\end{assumption}

\begin{assumption}
\label{as:density}
For every row $\va_i$ and every unit vector $\mathbf u\in\mathbb R^n$, the
random variable $\langle\va_i,\mathbf u\rangle$ has a density bounded by
$D\sqrt n$.
\end{assumption}

For example, rows drawn independently and uniformly from the unit sphere,
equivalently normalized Gaussian rows, satisfy both assumptions for $n\ge3$
with $K=\sqrt{8/3}$ and $D=1/\sqrt{2\pi}$ (as follows from~\cite[Section~2.5]{vershynin2018high}).

Let $C_K$ denote the constant of~\cite[Theorem~3.1]{haddock2022quantilebased},
which depends only on $K$, and let $C_\flat=(q-\beta)^{3/2}/(24D)$. For
$\gamma\in(0,1/2]$ and $\delta>1$, we define
\[
\begin{aligned}
\mu_*&=C_\flat^2-2\sqrt\beta\,\omega\,C_K^2,
&\qquad
c_\eta&=\frac{2(1-\gamma)(1-\delta^{-1})\,\mu_*}{q},\\
\alpha_*&=\frac{2\gamma\,\mu_*}{C_K^2\,\Lambda},
&\qquad
C_\xi&=\frac{8\delta}{1-\gamma}\Bigl(\frac{\sqrt\beta}{\omega\,\mu_*}
+\frac{\gamma}{2\omega^2}\Bigr),
\end{aligned}
\]
where $\Lambda=1+4\beta\omega^2+\bigl(C_K^2+4\beta\omega(\beta\omega+C_K)\bigr)/q$ and $\omega=1/\sqrt{1-q-\beta}$.

\begin{corollary}[Random design]
\label{cor:random_design}
Let Assumptions~\ref{as:model}--\ref{as:density} hold, and
suppose that $\mu_*>0$ and $m/n>C_0(q-\beta)^{-1}\log(DK/(q-\beta))$. Then,
with probability at least $1-2e^{-cm}-3e^{-c_qm}$ over the draw of $\mA$, the
following holds for all $\gamma\in(0,1/2]$, $\delta>1$, $s\ge1$ and
$0<\alpha\le\alpha_*$: if $C_{\rm corr}>2r_\infty\bigl(1+\sqrt{C_\xi\,n}\bigr)$,
then, for all $k\ge0$,
\[
\mathbb E\bigl[\|\ve_k\|^2\,\big|\,\mA\bigr]
\le\Bigl(1-\frac{c_\eta\,\alpha}{n}\Bigr)^k\|\ve_0\|^2
+\frac{\alpha\,n\,r_\infty^2}{c_\eta\,s}.
\]
Here $C_0$ and $c$ are positive absolute constants.
\end{corollary}

Thus the root-mean-square error horizon is of order $\sqrt{\alpha n/s}\,r_\infty$ and
the separation requirement of order $\sqrt n\,r_\infty$. For $s=1$ and $\alpha$ of order one,
the horizon has the same order as for QuantileRK with bounded
noise~\cite[Theorem~II.1]{jarman2021quantilerk}.

\begin{proof}
\proofstep{Step 1: spectral bounds.} Under Assumption~\ref{as:subgaus},
\cite[Theorem~3.1]{haddock2022quantilebased} gives
$\sigma_{\max}\le C_K\sqrt{m/n}$ with probability at least $1-2e^{-cm}$.
Under Assumption~\ref{as:density}, if
$m/n>C_0(q-\beta)^{-1}\log(DK/(q-\beta))$, the proof
of~\cite[Proposition~3.4]{haddock2022quantilebased} gives
$\sigma_{\min,q-\beta}\ge C_\flat\sqrt{m/n}$ with
$C_\flat=(q-\beta)^{3/2}/(24D)$, with probability at least $1-3e^{-c_qm}$. On the event that both these bounds hold, 
the restricted condition number satisfies
$\sigma_{\min,q-\beta}^2/\sigma_{\max}^2\ge C_\flat^2/C_K^2$, independently of
the dimensions.

\proofstep{Step 2: rate and stepsize.} By Step~1, with probability at least $1-2e^{-cm}-3e^{-c_qm}$, 
$\mu\ge\mu_*m/n$ and $\kappa\le C_K^2/(qn)$, so the definitions
in~\eqref{eq:main_constants} give $\eta\ge c_\eta/n$ and
$\alpha_{\max}(s)\ge\alpha_{\max}(1)\ge\alpha_*$ for every $s\ge1$.

\proofstep{Step 3: separation.} It remains to bound $\xi_s$ by a
dimension-free multiple of $\sqrt n\,r_\infty$. By~\eqref{eq:main_constants},
$\xi_s^2$ is proportional to
$(m/\mu)\bigl(\sqrt\beta+\beta\,\alpha_{\max}(s)\,(\omega/s
+(\beta\omega+\sqrt{q\kappa})/q)\bigr)$. Note that
\[
\alpha_{\max}(s)=\frac{2\gamma\mu/\sigma_{\max}^2}{\Delta(s)},
\quad
\Delta(s)=\frac{1+4\beta\omega^2}{s}+\kappa
+\frac{4\beta\omega(\beta\omega+\sqrt{q\kappa})}{q} \ge 4\omega\beta\bigl(\frac{\omega}{s}+\frac{\beta\omega+\sqrt{q\kappa}}{q}\bigr),
\]
and therefore
\[
\beta\,\alpha_{\max}(s)\Bigl(\frac\omega s
+\frac{\beta\omega+\sqrt{q\kappa}}{q}\Bigr)
\le\frac{\gamma\mu}{2\omega\,\sigma_{\max}^2}.
\]
Inserting this bound into the definition of $\xi_s$, and using
$m/\mu\le n/\mu_*$ and, for unit rows, $\sigma_{\max}^2\ge\|\mA\|_F^2/n=m/n$,
gives
\[
\xi_s^2\le\frac{8\delta\,r_\infty^2}{1-\gamma}
\Bigl(\frac{\sqrt\beta}{\omega}\,\frac m\mu
+\frac{\gamma}{2\omega^2}\,\frac{m}{\sigma_{\max}^2}\Bigr)
\le C_\xi\,n\,r_\infty^2 .
\]
Therefore $C_{\rm corr}>2r_\infty(1+\sqrt{C_\xi n})$ implies the separation
condition~\eqref{eq:separation} for every $s\ge1$.

\end{proof}

\subsection{Further beyond the horizon}
\label{subsec:schedules}
Throughout this subsection we write $u_k=\mathbb E\|\ve_k\|^2$, where the
expectation is over all randomness up to iteration $k$; $u_0=\|\ve_0\|^2$.

Theorem~\ref{thm:main} breaks the error horizon of QuantileRK: AQRK's mean-square horizon is
proportional to $\alpha/s$ and can be made arbitrarily small by the choice of
$\alpha$ and $s$; however, for fixed parameters it remains positive. Growing $s$ or decreasing $\alpha$ can make the horizon vanish:
\begin{corollary}[Vanishing horizon]
\label{cor:vanishing}
Let the assumptions of Theorem~\ref{thm:main} hold for some
$\gamma\in(0,1/2]$, $\delta>1$ and a baseline batch size $s_0\ge1$, including
the separation condition $C_{\rm corr}>2(r_\infty+\xi_{s_0})$. Run
Algorithm~\ref{alg:framework} with deterministic stepsizes
$0<\alpha_k\le\alpha_{\max}(s_0)$ and deterministic integer batch sizes
$s_k\ge s_0$. If
\[
\sum_{k=0}^{\infty}\alpha_k=\infty
\qquad\text{and}\qquad
\frac{\alpha_k}{s_k}\to0,
\]
then $u_k\to0$. In particular, this holds for a fixed batch with stepsizes
$\alpha_k\to0$ and $\sum_k\alpha_k=\infty$, and for a fixed stepsize with
$s_k\to\infty$.
\end{corollary}
\begin{proof}
By Proposition~\ref{prop:onestep}, after taking total expectations,
\[
u_{k+1}\le(1-q_k)\,u_k+q_kh_k,
\qquad
q_k=\eta\alpha_k,
\qquad
h_k=\frac{\alpha_k r_\infty^2}{\eta s_k},
\]
where $0<q_k<1$ as in the proof of Theorem~\ref{thm:main}. Let
$\varepsilon'>0$. Since $h_k\to0$, there is $J$ with $h_k\le\varepsilon'$ for
all $k\ge J$, and then
$(u_{k+1}-\varepsilon')_+\le(1-q_k)(u_k-\varepsilon')_+$. Hence
$(u_k-\varepsilon')_+\le(u_J-\varepsilon')_+\prod_{j=J}^{k-1}(1-q_j)$ for
$k>J$, and the product tends to zero because $\sum_jq_j=\infty$. Thus
$\limsup_k u_k\le\varepsilon'$, and since $\varepsilon'$ is arbitrary,
$u_k\to0$.
\end{proof}

The two standard ways of satisfying these conditions also come with explicit
rates.

\paragraph{Decreasing stepsizes.} For a fixed batch, the schedule, adapted from~\cite{gower2019sgd}, results in an initial phase of exponential decay, at the rate of
Theorem~\ref{thm:main}, followed by an $O(1/\sqrt{k})$ convergence to zero phase; the switch $k^*$ is the first
index when the decreasing stepsize does not exceed $\alpha_0$:
\begin{corollary}[Decreasing stepsizes]
\label{cor:diminishing_stepsize}
Let the assumptions of Theorem~\ref{thm:main} hold for some
$\gamma\in(0,1/2]$, $\delta>1$ and $s\ge1$, including the separation
condition $C_{\rm corr}>2(r_\infty+\xi_s)$, and let
$0<\alpha_0\le\alpha_{\max}(s)$. Let $k^*$ be the smallest integer $k\ge0$
such that
\[
\frac{2k+1}{(k+1)^2}\le\eta\alpha_0,
\]
and run Algorithm~\ref{alg:framework} with $s_k\equiv s$ and
\[
\alpha_k=
\begin{cases}
\alpha_0, & k<k^*,\\[1ex]
\dfrac{2k+1}{\eta(k+1)^2}, & k\ge k^*.
\end{cases}
\]
Then $k^*\ge1$ and, for every $k\ge k^*$,
\[
u_k\le\frac{(k^*)^2}{k^2}\Bigl((1-\eta\alpha_0)^{k^*}u_0
+\frac{\alpha_0r_\infty^2}{\eta s}\Bigr)
+\frac{4(k-k^*)\,r_\infty^2}{\eta^2s\,k^2}.
\]
In particular, $u_k\le(k^*/k)^2u_{k^*}+4r_\infty^2/(\eta^2sk)=O(1/k)$.
\end{corollary}

Note that averaging batch size $s$ helps in both phases; in particular, the constant of the squared error decay 
$O(1/k)$ term is inversely proportional to the batch size $s$.
\begin{proof}
\emph{Admissible stepsizes.} The ratio $(2k+1)/(k+1)^2$ equals $1$ at $k=0$
and decreases with $k$. Since $\eta\alpha_0<1$, this gives $k^*\ge1$ and
$\alpha_k\le\alpha_0$ for $k\ge k^*$, so $\alpha_k\le\alpha_0\le\alpha_{\max}(s)$
for all $k$. So, Proposition~\ref{prop:onestep} applies to give
\begin{equation}
\label{eq:unif-recursion}
u_{k+1}\le(1-\eta\alpha_k)\,u_k+\frac{\alpha_k^2r_\infty^2}{s}
\qquad\text{for all }k\ge0.
\end{equation}

\emph{Initial phase.}  For $k<k^*$, the stepsize is $\alpha_0$, and iterating
\eqref{eq:unif-recursion} gives
\begin{equation}\label{eq:initial-phase}
u_{k^*}\le(1-\eta\alpha_0)^{k^*}u_0
+\frac{\alpha_0r_\infty^2}{\eta s}.
\end{equation}

\emph{Decreasing phase.} For $k\ge k^*$, the choice of $\alpha_k$ gives
\[
1-\eta\alpha_k=\frac{k^2}{(k+1)^2},
\qquad
\alpha_k^2=\frac{(2k+1)^2}{\eta^2(k+1)^4}\le\frac{4}{\eta^2(k+1)^2},
\]
where we used $2k+1\le2(k+1)$. Inserting these
into~\eqref{eq:unif-recursion} and multiplying by $(k+1)^2$ gives
$(k+1)^2u_{k+1}\le k^2u_k+4r_\infty^2/(\eta^2s)$. Summing over
$k=k^*,\dots,K-1$ for some $K>k^*$, the left- and right-hand sides telescope,
and $K^2u_K\le(k^*)^2u_{k^*}+4(K-k^*)\,r_\infty^2/(\eta^2s)$, which together
with~\eqref{eq:initial-phase} proves Corollary~\ref{cor:diminishing_stepsize}.
\end{proof}

\paragraph{Increasing batches.}
Alternatively, the stepsize can stay fixed while the number of fresh update
observations grows, as known in the literature for data fitting and
stochastic optimization~\cite{friedlander2012hybrid,smith2018dont}. Here, this
works after quantile filtering, despite the persistent sparse corruption. Any
$s_k\to\infty$ gives $u_k\to0$ by Corollary~\ref{cor:vanishing}, and geometric
growth gives a geometric rate:

\begin{corollary}[Increasing update batches]
\label{cor:increasing_batch}
Let the assumptions of Theorem~\ref{thm:main} hold for some
$\gamma\in(0,1/2]$, $\delta>1$ and a baseline batch size $s_0\ge1$,
including the separation condition $C_{\rm corr}>2(r_\infty+\xi_{s_0})$, and
fix $0<\alpha\le\alpha_{\max}(s_0)$. Run Algorithm~\ref{alg:framework} with
$\alpha_k\equiv\alpha$ and deterministic integer batch sizes $s_k\ge s_0$, and
write $\theta=1-\eta\alpha$ and $\nu=\alpha^2r_\infty^2.$
Then $\theta\in(0,1)$ and
\begin{equation}
\label{eq:batch_recursion}
u_k\le\theta^ku_0+\nu\sum_{j=0}^{k-1}\frac{\theta^{k-1-j}}{s_j}.
\end{equation}
For the geometric schedule
$s_k=\lceil s_0\zeta^{-k}\rceil$ with $\theta<\zeta<1$,
\begin{equation}
\label{eq:geometric_batch_bound}
u_k\le\theta^ku_0+\frac{\nu}{s_0}\,\frac{\zeta^k-\theta^k}{\zeta-\theta}
\le\Gamma\zeta^k,
\qquad
\Gamma=u_0+\frac{\nu}{s_0(\zeta-\theta)}.
\end{equation}
\end{corollary}

\begin{proof} By Proposition~\ref{prop:onestep}, after taking total expectations,
$u_{k+1}\le\theta u_k+\nu/s_k$, and $0<\eta\alpha<1$ as in the proof of
Theorem~\ref{thm:main}, so $\theta\in(0,1)$. Induction on $k$
gives~\eqref{eq:batch_recursion}. For $s_j=\lceil s_0\zeta^{-j}\rceil$, we
have $1/s_j\le\zeta^j/s_0$ and
\[
\sum_{j=0}^{k-1}\theta^{k-1-j}\zeta^j
=\zeta^{k-1}\sum_{j=0}^{k-1}\Bigl(\frac{\theta}{\zeta}\Bigr)^{k-1-j}
=\frac{\zeta^k-\theta^k}{\zeta-\theta}.
\]
Inserting this into~\eqref{eq:batch_recursion} gives the first inequality
in~\eqref{eq:geometric_batch_bound}; the second follows from
$\theta^k\le\zeta^k$ and $\zeta^k-\theta^k\le\zeta^k$.
\end{proof}

\paragraph{Smaller steps or larger batches?}
\label{subsec:cost}

While both approaches drive the error to zero, decreasing the stepsize decreases the convergence rate, and increasing the batch sizes makes every iteration more expensive. By Corollaries~\ref{cor:diminishing_stepsize}
and~\ref{cor:increasing_batch}, both schedules reach $u_K =\mathbb E\|\ve_K\|^2\le\varepsilon^2$
with $O(\varepsilon^{-2})$ fresh update observations: decreasing stepsizes with
a fixed batch $s$ after $O(\varepsilon^{-2}/s)$ iterations, and geometric
batches after $K=O(\log(1/\varepsilon))$ iterations, since then
$\sum_{k<K}s_k=O(\zeta^{-K})=O(\varepsilon^{-2})$ by~\eqref{eq:geometric_batch_bound}. However, the exact comparison depends on the specifics of implementation: if the update observations reuse the cached inner
products from the screening stage of the same iteration (Remark~\ref{rem:sampling_stages}), fresh
observations of the noisy right-hand side may be inexpensive relative to screening (e.g., if every observation costs $O(1)$, then a batch of size $O(m_0n)$ has the complexity of the order of the corresponding initial screening stage with $m_0$ rows, $m_0 \le m$), which makes increasing batches especially attractive. 

In practice, there is little point in decreasing the stepsize or enlarging the batch from the start. Initially, the convergence behavior is dominated by the initial error rather than by noise, so doing either early would lead to unnecessary compute expense. A reasonable policy can be to keep the
parameters constant initially, then increase the batch, and decrease the stepsize if further batch growth becomes too expensive.

This schedule theory assumes that the noise is resampled at every observation. If the noise cannot be resampled, neither schedule provably removes the effect of noise after the bias introduced by the filtering stage. In
practice, noise can combine a part that changes between observations with
a part that does not, such as calibration offsets or model mismatch, and the
latter acts like fixed noise. Our experiments below test both resampled
noise, where Theorem~\ref{thm:main} applies, and fixed noise, which lies
beyond it.

\section{Experiments}
\label{sec:experiments}

First, in Section~\ref{sec:general_conv}, we demonstrate the conclusions of Theorem~\ref{thm:main} on synthetic data, also testing the regimes beyond the assumptions of the theory. Then, we demonstrate the applicability and advantages of AQRK on a real-world tomography dataset (Section~\ref{sec:tomography_exp}) and as a robust solver for alternating least-squares inner problems  inside a low-rank tensor fitting task (Section~\ref{subsec:cp_application}).

Throughout the section, we compare
\emph{AQRK (fresh)} and \emph{AQRK (reused)}, the two variants of
Algorithm~\ref{alg:framework}, and \emph{AQRK-all}, the update~\eqref{eq:aqrk-all},
with three baselines: randomized Kaczmarz (\emph{RK}); averaged randomized Kaczmarz
(\emph{ARK}), which averages $s$ rows drawn uniformly from $[m]$ with
replacement; and QuantileRK (\emph{QRK}) implemented as the
reused variant of Algorithm~\ref{alg:framework} with $s=1$ for fair comparison. All quantile-based methods use full screening and the stepsize $\alpha=1$ unless otherwise specified.
Each experiment uses $10$ trials, with the same systems for all methods within
a trial. We report the relative error $\|\vx_k-\vx^\star\|/\|\vx^\star\|$ from
$\vx_0=\boldsymbol 0$. Every figure shows the pointwise median over the trials,
with a band from the $25$th to the $75$th percentile  at each iteration.

All code used to produce the results in this section is publicly available at
\url{https://github.com/shvaiko-s/aqrk} (version \texttt{49f80bc}), including the
implementations of all methods, the scripts used to run the experiments and to
generate the figures and tables, and the data and random seeds needed to reproduce
them. The experiments were run in Python~$3.9$ on a laptop with an Apple M4 processor
running macOS~$15$.

\subsection{Synthetic systems}
\label{sec:general_conv}
In this section, each row of $\mA\in\mathbb R^{2000\times100}$ is a
standard Gaussian vector normalized to unit length,
$\vx^\star\sim N(\boldsymbol 0,\mathbf I_n)$ and $\vb^\star=\mA\vx^\star$. A fraction $\beta=0.3$ of the
equations, that is, $600$ rows, is corrupted, and we take $q = 0.6$. Unless stated otherwise, the
query noise is drawn independently from $\mathrm{Uniform}[-0.02,0.02]$.

\paragraph{Corruption models.}
\begin{itemize}
\item \emph{Oblivious.} The support is drawn uniformly without replacement, and
$\rho_i\sim\mathrm{Uniform}[a,b]$ independently.
\item \emph{Sign-flip.} The support is drawn uniformly without replacement, and
$\rho_i=-2b_i^\star$. These corruptions might not satisfy the separation condition.
\item \emph{Singular-vector-based.} Let $\mathbf u$ be the left singular vector of
$\mA$ for its smallest singular value. The support consists of the $600$
indices with the largest $|u_i|$, and $\rho_i=f\,\mathrm{sign}(u_i)$ for a
fixed magnitude $f$. 
\end{itemize}

\subsubsection{Base experiment: fresh versus reused updates}
\label{sec:base_exp}
\begin{figure}[htbp]
    \centering
    \includegraphics[width=\linewidth]{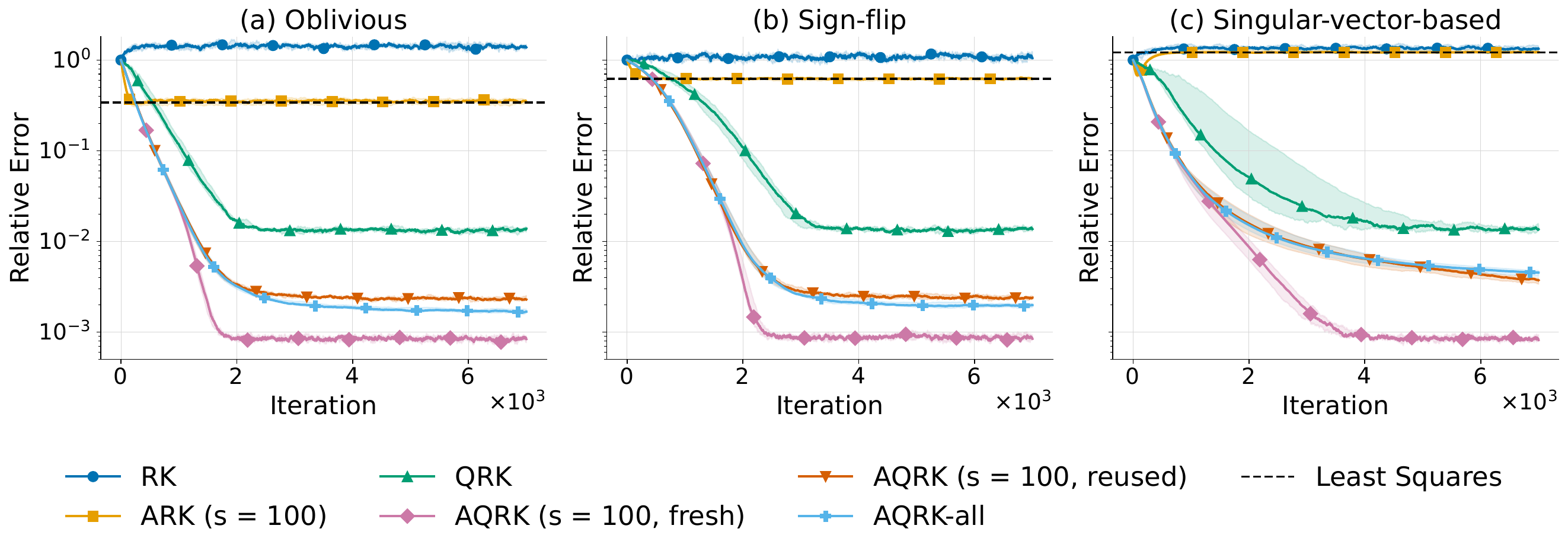}
 \caption{Relative error with query noise under (a) oblivious uniform
corruption on $[2,3]$, (b) sign-flip corruption, and (c) singular-vector
corruption with $f=2$. The reference line is the median relative error
of the least-squares solution for the corrupted right-hand side without
query noise, the unfiltered baseline}
    \label{fig:generic_convergence}
\end{figure}

In Fig.~\ref{fig:generic_convergence}, we show how re-querying the selected
equations removes the filtering bias. We can see that for harder corruption models the convergence can be slower but the horizon values remain consistent.
In all three corruption types, Fig.~\ref{fig:generic_convergence}(a)--(c),
RK and ARK, which do not filter corruptions, stall, every
AQRK variant levels off below QRK, and
\emph{AQRK ($s=100$, fresh)} attains the lowest floor. \emph{AQRK-all} averages all $h=1200$ eligible rows, yet it gains little compared to the reused sample of only $100$ rows.
Averaging more rows in the re-queried regime meaningfully decreases the horizon as we see in the next figure.

\subsubsection{Choosing update batch size and step size}
\label{sec:batch_exp}
\begin{figure}[htbp]
    \centering
    \begin{minipage}[t]{0.48\linewidth}
        \centering
        \includegraphics[width=\linewidth]{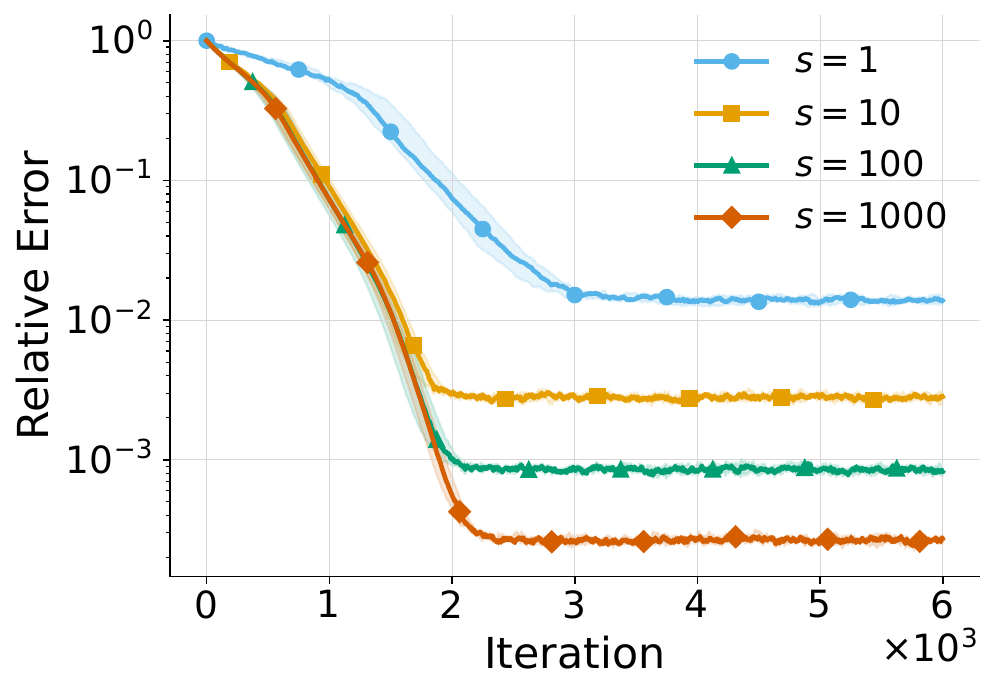}
        \caption{Effect of the update batch size $s$ on the horizon of \emph{AQRK (fresh)}}
        \label{fig:s_sweep}
    \end{minipage}\hfill
    \begin{minipage}[t]{0.48\linewidth}
        \centering
        \includegraphics[width=\linewidth]{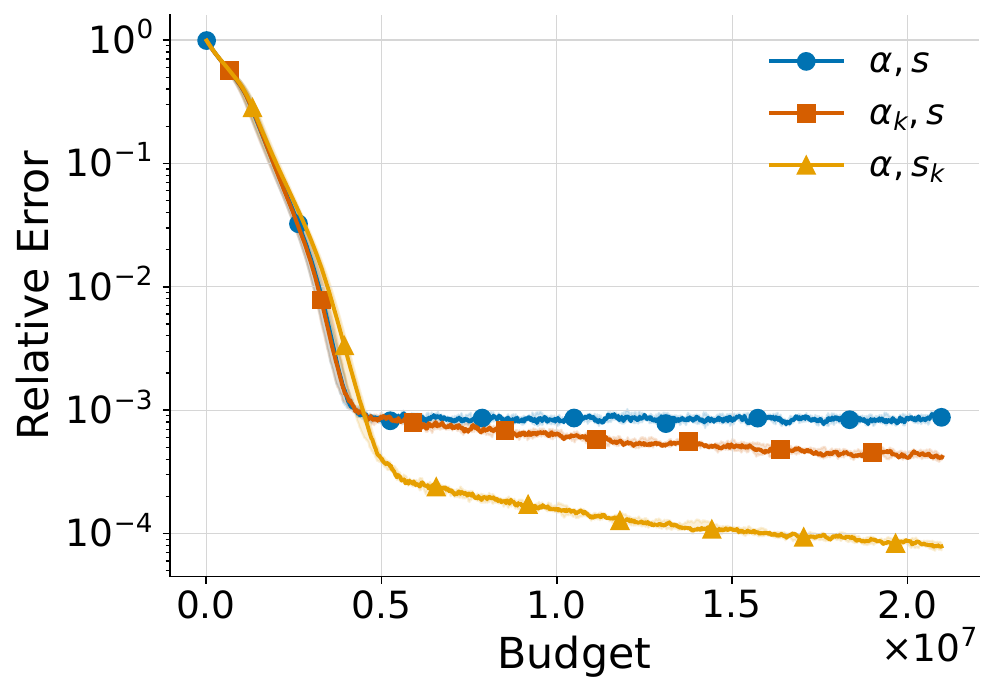}
        \caption{\emph{AQRK (fresh)} with constant and scheduled stepsizes and batch sizes: here, $\alpha=1$, $s=100$;
        $\alpha_k=k^*/k$ for $k\ge k^*=2500$ and $\alpha_k = 1$ otherwise; $s_k=\lceil100\cdot0.999^{-k}\rceil$}
        \label{fig:schedules}
    \end{minipage}
\end{figure}

From Theorem~\ref{thm:main}, the root-mean-square
error horizon should decrease as $O(s^{-1/2})$. We test this prediction with
\emph{AQRK (fresh)} for $s\in\{1,10,100,1000\}$ under oblivious corruption
on $[1,1.2]$. Indeed, Fig.~\ref{fig:s_sweep} shows the relative-error floors at about
$1.4\times10^{-2}$, $2.8\times10^{-3}$, $8.5\times10^{-4}$ and
$2.7\times10^{-4}$ for $s=1,10,100,1000$. Each tenfold increase of $s$
beyond $s=10$ lowers the floor by a factor of about $3.2$, close to the
predicted $\sqrt{10}\approx3.16$. 
Additionally, Fig.~\ref{fig:s_sweep} shows that for $s\ge10$ the curves nearly coincide during the initial
descent and level off faster than the curve of $s=1$. 

Corollaries~\ref{cor:diminishing_stepsize} and~\ref{cor:increasing_batch}
propose to remove the error horizon by decreasing the stepsize or by increasing the
batch size. In Fig.~\ref{fig:schedules}, we compare a form of the schedule from
Corollary~\ref{cor:diminishing_stepsize} with $k^* = 2500$, a form of the schedule from
Corollary~\ref{cor:increasing_batch} with $\zeta=0.999$ and a constant schedule, using \emph{AQRK (fresh)} with
oblivious corruption on $[1,1.2]$. We measure cost in
residual queries: iteration $k$ queries all $m=2000$ rows for screening and
$s_k$ rows for the update. Each schedule runs until it would exceed a total
of $2.1\times10^7$ queries, which corresponds to $10^4$ iterations with
$s=100$. Fig.~\ref{fig:schedules} shows the error against the number of residual
queries, $\sum_k(m+s_k)$.
Counting queries treats screening and update queries alike. 
Consistent with Section~\ref{subsec:cost}, increasing the batch size reduces
the error faster per query than decreasing the stepsize, achieving an error lower by a factor of approximately $5.2$.

\subsubsection{Sensitivity to corruption magnitude}
\label{sec:magnitude_exp}

\begin{figure}[htbp]
    \centering
        \includegraphics[width=\linewidth]{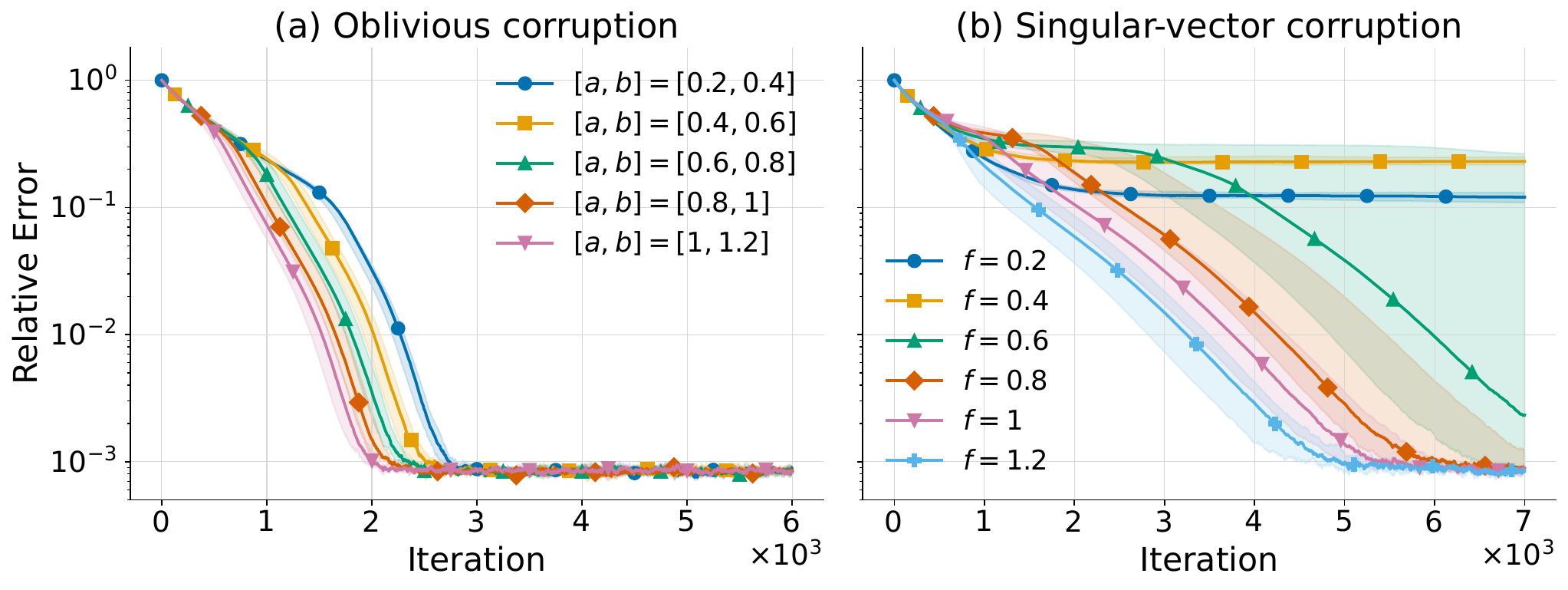}
    \caption{Relative error of \emph{AQRK (fresh)} under varying corruption magnitude for oblivious (a) and singular-vector-based (b) corruption models}
    \label{fig:magnitude_sweep}
\end{figure}

Small corruptions are harder to filter, because their residuals can be
comparable to those of uncorrupted rows.
To investigate robustness as the corruption magnitude and separation decreases, we vary the
corruption magnitude over a broad range while keeping the corruption fraction
and noise model fixed.
We test empirically necessary separation on two corruption models: oblivious and adversarial
singular-vector-based corruption, see Fig.~\ref{fig:magnitude_sweep}. Under oblivious corruption (Fig.~\ref{fig:magnitude_sweep}(a)), every magnitude range reaches the same error floor just below $10^{-3}$, and smaller magnitudes only delay convergence slightly.
Singular-vector corruption is harder (Fig.~\ref{fig:magnitude_sweep}(b)). For
$f\le2\sqrt n\,r_\infty$ the corrupted rows are not filtered out, and the
error stagnates at a level that grows with $f$. For larger $f$ the iterates
can leave this plateau and reach the same floor as under oblivious
corruption, but the time spent on the plateau grows, and varies more
between trials, as $f$ decreases.  Both results
are consistent with the order $\sqrt n\,r_\infty$ in
Corollary~\ref{cor:random_design}, but the constant depends on how the
corruption is aligned with $\mA$.

\subsection{Tomography test problem}
\label{sec:tomography_exp}

To test the methods on a structured, ill-conditioned matrix coming from real-world data, we use the
tomography test problem of Regularization
Tools~\cite{hansen2007regularization}, with the image provided by the test
problem as the reference target $\vx^\star$. After normalizing the rows and
removing zero rows, the system has $m\times n=1148\times100$ and condition
number about $62$, and we set $\vb^\star=\mA\vx^\star$, whose entries are at
most about $1$ in magnitude. In each trial, $137$ entries, a fraction
$\beta\approx0.12$ (and we choose $q\approx0.78$ for the experiments), are corrupted by $\rho_i=\pm f_i$ with independent random
signs and $f_i\sim\mathrm{Uniform}[5,6]$. The noise is Gaussian with standard
deviation $\sigma=0.06$, in two settings: with query noise, every
residual query receives an independent draw; with fixed noise, a
single draw is added to $\vb$ and the queries are noise-free. We note that the underlying system is significantly worse conditioned than a Gaussian system of that size, and the Gaussian noise
is chosen as a natural modeling assumption, but it is unbounded, so Assumption~\ref{as:model} does not hold.

\begin{figure}[htbp]
    \centering
    \includegraphics[width=\linewidth]{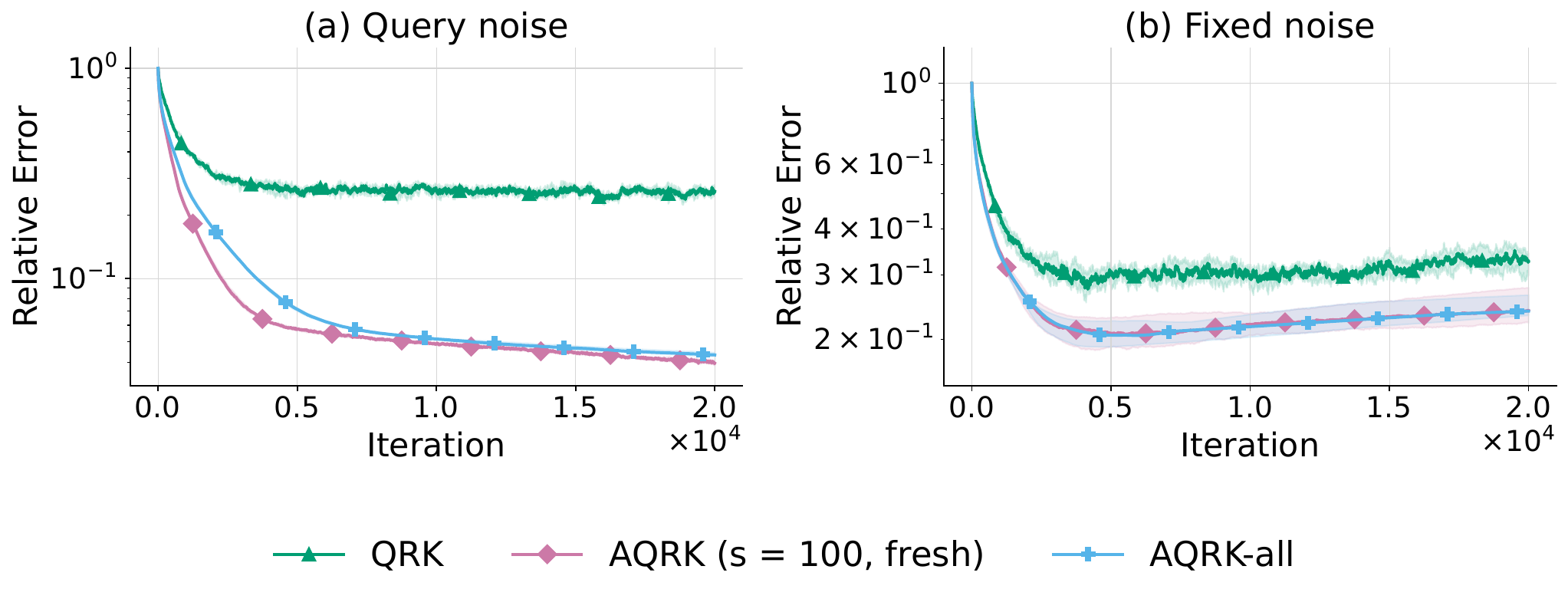}
    \caption{Relative error for the tomography problem under (a) Gaussian query noise and
    (b) fixed Gaussian noise. RK and ARK
    are not shown; their relative errors exceed~$1$}
    \label{fig:tomography}
\end{figure}

With query noise (Fig.~\ref{fig:tomography}(a)), QRK levels off at a relative
error of about $0.25$, whereas \emph{AQRK ($s = 100$, fresh)} and \emph{AQRK-all} reach about $0.04$ after
$20000$ iterations and are still decreasing. On this ill-conditioned matrix,
the error at the end of the run is limited by slow convergence rather than by
the error floor, and the difference between \emph{AQRK (fresh)} and \emph{AQRK-all} is small. With
fixed noise (Fig.~\ref{fig:tomography}(b)), the averaging also helps reduce the horizon compared to the QRK results; however, all the methods show
semi-convergence, as the error reaches a minimum after about $5000$ iterations
and then increases towards that of the least-squares solution of the
uncorrupted equations. RK and
ARK, which do not filter and diverge, are not shown.

In summary, quantile filtering remains essential and
averaging over resampled observations improves the attainable accuracy of
QRK substantially on this structured matrix, even though the data fails to satisfy most of the conditions of Theorem~\ref{thm:main}.

\subsection{Filtered inner iterations for CP tensor fitting}
\label{subsec:cp_application}

In this section, the problem of CP tensor fitting tests AQRK as a component of a nonlinear numerical
computation. A rank-$R$ CP representation has the form
$\mathcal X=\sum_{r=1}^R{\bf a}_r\circ{\bf b}_r\circ{\bf c}_r
=\llbracket{\bf A},{\bf B},{\bf C}\rrbracket$,
where $\circ$ denotes the vector outer product~\cite{kolda2009tensor}. CP alternating least squares (CP-ALS) cyclically updates one factor
while the other two are fixed. With ${\bf F}_1={\bf A}$,
${\bf F}_2={\bf B}$, ${\bf F}_3={\bf C}$, the mode-$j$ unfolding
${\bf Y}_{(j)}$ of the observed tensor $\mathcal Y$ is fitted as
\begin{equation}\label{eq:cp_unfold}
{\bf Y}_{(j)}\approx{\bf F}_j{\bf M}_j^T,\qquad
{\bf M}_1={\bf C}\odot{\bf B},\quad
{\bf M}_2={\bf C}\odot{\bf A},\quad
{\bf M}_3={\bf B}\odot{\bf A},
\end{equation}
where $\odot$ is the Khatri--Rao product. Each factor update therefore
splits into independent least-squares problems, one per row:
$({\bf F}_j)_{i,:}=\vx^T$ solves
${\bf M}_j\vx\approx({\bf Y}_{(j)})_{i,:}^{\top}$.
Algorithm~\ref{alg:cp-fitting} replaces these solves by filtered,
averaged updates, warm-started from the current factors, with
coefficient rows and right-hand sides normalized once per factor update.
\begin{algorithm}[htbp]
\caption{CP fitting with filtered, averaged inner updates}
\label{alg:cp-fitting}
\begin{algorithmic}[1]
\Require Observed tensor $\mathcal Y$, initial factors
         ${\bf F}_1,{\bf F}_2,{\bf F}_3$, outer sweeps $N_{\rm out}$
\Require AQRK parameters, fresh/reused update choice,
         and inner stopping rule
\For{$\ell=1,\ldots,N_{\rm out}$}
    \For{$j=1,2,3$}
        \State Form ${\bf M}_j$ as in~\eqref{eq:cp_unfold}; drop its zero
               rows and the corresponding columns of ${\bf Y}_{(j)}$
        \State ${\bf D}_j\gets\operatorname{diag}\bigl(\|({\bf M}_j)_{t,:}\|_2\bigr)$,
               $\widetilde{\bf M}_j\gets{\bf D}_j^{-1}{\bf M}_j$,
               $\widetilde{\bf Y}_{(j)}\gets{\bf Y}_{(j)}{\bf D}_j^{-1}$
               \Comment{Unit-norm rows}
        \For{each row $i$ of ${\bf F}_j$}
            \State $({\bf F}_j)_{i,:}\gets\operatorname{AQRK}
                   \bigl(\widetilde{\bf M}_j,\,
                   (\widetilde{\bf Y}_{(j)})_{i,:}^{\top},\,
                   ({\bf F}_j)_{i,:}^{\top}\bigr)^{\top}$
                   \Comment{Warm start}
        \EndFor
    \EndFor
    \State Normalize the columns of ${\bf F}_2,{\bf F}_3$, rescaling
           those of ${\bf F}_1$ to preserve the tensor
\EndFor
\State \Return $\llbracket{\bf F}_1,{\bf F}_2,{\bf F}_3\rrbracket$
\end{algorithmic}
\end{algorithm}

We demonstrate Algorithm~\ref{alg:cp-fitting} on synthetic 
$30\times30\times30$ tensors. The exact target tensor
$\mathcal X_{\rm ex}$ is modeled with independent standard Gaussian factors and $R = 5$ components. To introduce sparse biased corruptions, we replace $10\%$ of $\mathcal X_{\rm ex}$ entries, chosen uniformly without replacement, with the new entries sampled from the $\mathrm{Uniform}[-7.3,7.3]$ distribution. We note that these corrupted values cannot be identified simply based on their size as the original tensor's entry range is
approximately $[-32,18]$.

\begin{figure}[htbp]
\centering
\includegraphics[width=\linewidth]{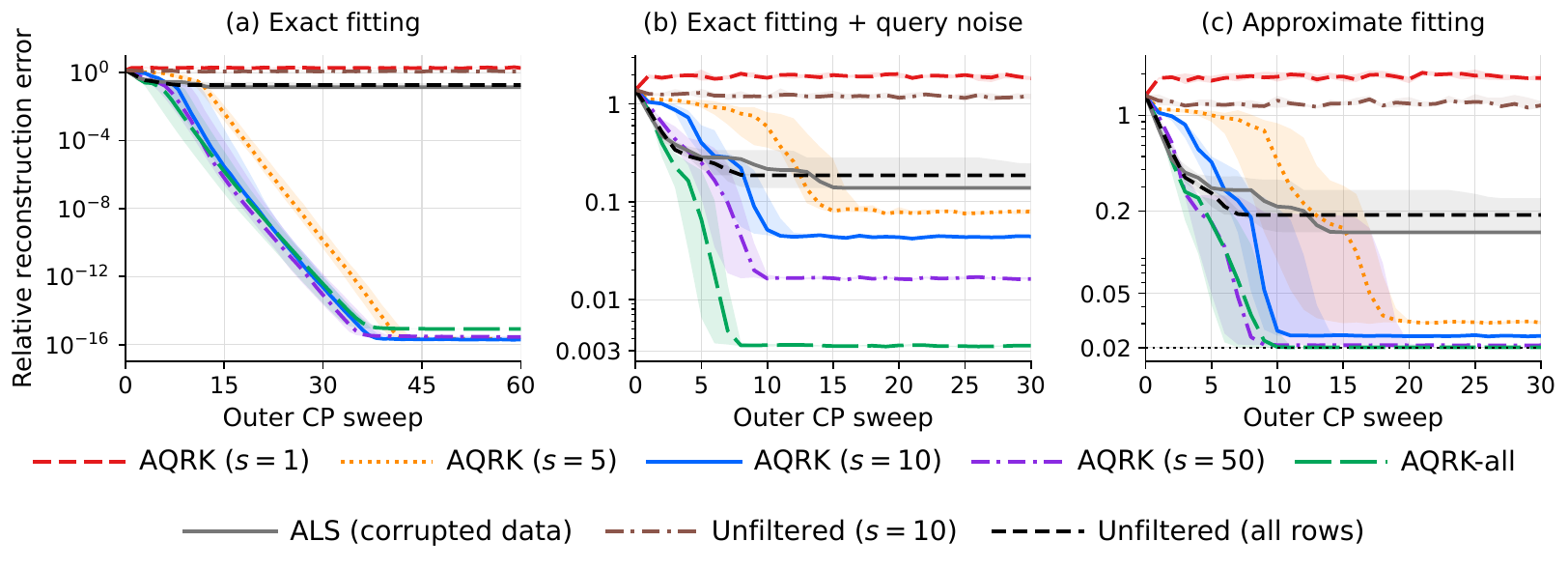}
\caption{CP fitting with $10\%$ corrupted entries: (a) exact target,
(b) exact target with query noise, and (c) approximately low-rank target. The relative reconstruction error
is $
\|\widehat{\mathcal X}-\mathcal X_\star\|_F/\|\mathcal X_\star\|_F$,
where $\widehat{\mathcal X}$ is the current CP reconstruction and
$\mathcal X_\star$ is the uncorrupted target: $\mathcal X_{\rm ex}$ in
Fig.~\ref{fig:cp_experiment}(a),(b) and $\mathcal X_{\rm app}$ in
Fig.~\ref{fig:cp_experiment}(c). The dotted line at $0.02$ in (c) is the error
$\|\mathcal X_{\rm app}-\mathcal X_{\rm ex}\|_F/\|\mathcal X_{\rm app}\|_F$
of the generating rank-$5$ tensor}
\label{fig:cp_experiment}
\end{figure}

First, we aim to recover $\mathcal X_{\rm ex}$ from this corrupted version. Note that even for the exact low-rank target, the clean inner systems are inconsistent while the other factors are inaccurate and change between outer sweeps, so already the first setting lies beyond our linear-system theory. Figure~\ref{fig:cp_experiment}(a) shows that AQRK with $s=10$ recovers the exact target in all ten trials, whereas QRK
recovers it in none. Unfiltered averaging, even over all rows, removes
row-sampling fluctuations but not the fixed corruption, and stalls without convergence.

Second, we consider the case when every measurement also includes a  query noise, modeled as additional independent Gaussian perturbations with standard deviation
$\sigma=0.05\|\mathcal X_{\rm ex}\|_F/\sqrt{30^3}$ on every entry query.
Figure~\ref{fig:cp_experiment}(b) shows that
the errors level off at floors that decrease with the batch size, as
expected, and that averaging without filtering as well as filtered but not averaged iterations stall and ultimately give higher than $0.1$ reconstruction error.

Finally, in Fig.~\ref{fig:cp_experiment}(c), we consider a natural case when the underlying tensor is approximately rather than exactly low-rank. We define a new target as
\[
\mathcal X_{\rm app}
=\mathcal X_{\rm ex}
+0.02\,\frac{\|\mathcal X_{\rm ex}\|_F}{\|\mathcal G\|_F}\mathcal G,
\]
where $\mathcal G$ has independent standard Gaussian entries.
The perturbation is drawn once and held fixed across methods and trials;
replacement corruption is then applied to $\mathcal X_{\rm app}$. Here, note that $0.02$ is naturally the best expected reconstruction error as CP-ALS constructs an exactly $5$-component approximation. Algorithm~\ref{alg:cp-fitting} with  $s = 50$ essentially recovers a low-rank approximation of this accuracy (reaching $0.0208$ error after about $10$ outer iterations), the smaller block AQRK methods trail closely (requiring more outer iterations to reach a slightly worse error horizon), while unfiltered or not averaged methods do not decrease the error better than the standard ALS (which only reached $0.14$ median reconstruction error in our ten-trial experiment).

In each of ten trials, all methods start from the same initial factors,
which are also shared across the three cases. We use $q=0.85$, and re-query AQRK only for the experiment with query noise. The reference method ALS solves the least-squares systems with diagonal regularization
$10^{-8}{\bf I}_R$, while Kaczmarz-based solvers do not solve normal equations explicitly and so do not need regularization to avoid numerical singularities. Each inner iteration is stopped using the root-mean-square residual
$r_{\rm rms}$ of its system. Starting after $40$ updates, we evaluate
$r_{\rm rms}$ every $20$ updates and stop after two consecutive checks
with relative decrease below $0.1\%$, or once $r_{\rm rms}$ is at most
$10^{-12}$ times its first recorded value. Each inner iteration
performs at most $500$ updates.

\section{Conclusions and future directions}
\label{sec:conclusions}

Using AQRK (Algorithm~\ref{alg:framework}), we study how averaging and relaxation affect the attainable accuracy of
quantile-based solvers in the presence of sparse corruption and dense noise. The sparse
corruption may be non-centered and arbitrarily large above the separation
threshold. 
For independent repeated observations, we prove that under residual separation and
fresh post-selection queries, the AQRK algorithm yields an error horizon proportional to
$\sqrt{\alpha/s}$, as for averaging without corruption. In particular, for Gaussian-like data, with $s = O(n)$ and step size $1$, the horizon is of order the size of the noise $r_\infty$. Decreasing stepsizes or increasing independent batches give
mean-square recovery of the clean solution; in the repeated-observation model,
this answers the open question of~\cite{coria2024quantile}. Both schedules use
$O(\varepsilon^{-2})$ fresh observations to reach accuracy $\varepsilon$, but
increasing batches screen all rows only $O(\log(1/\varepsilon))$ times. The numerical experiments support the theory and additionally show that averaged quantile updates can
improve the attainable accuracy for more general errors of mixed scales, including fixed noise that cannot be resampled. 
AQRK also enables successful low-rank CP tensor fitting in cases where
alternating least-squares solvers, averaging techniques without filtering, as well as single-row robust updates  fail.

Future directions include further analysis of the fixed noise and low separation cases, as well as exploiting sparsity of the desired solution when it is applicable. It is interesting how the semi-convergence behavior observed in Fig.~\ref{fig:tomography}(b) depends on the parameters of the algorithm such as the batch and step size and how to pick the appropriate stopping rule. Optimal step- and batch-size schedules present another interesting direction.

Finally, for changing least-squares subproblems,
natural questions include stopping rules for the inner and outer loops and
reusing information accumulated in earlier outer iterations.

\backmatter

\bibliography{ref}
\end{document}